
%



\def\LaTeX{%
  \let\Begin\begin
  \let\End\end
  \def\Bcenter{\Begin{center}}
  \def\Ecenter{\End{center}}
  \let\Label\label
  \let\salta\relax
  \let\finqui\relax
  \let\futuro\relax}

\def\UK{\def\our{our}\let\sz s}
\def\USA{\def\our{or}\let\sz z}



\LaTeX

\USA


\salta

\documentclass[twoside,12pt]{article}
\setlength{\textheight}{24cm}
\setlength{\textwidth}{16cm}
\setlength{\oddsidemargin}{2mm}
\setlength{\evensidemargin}{2mm}
\setlength{\topmargin}{-15mm}
\parskip2mm


\usepackage{color}
\usepackage{amsmath}
\usepackage{amsthm}
\usepackage{amssymb}

\usepackage{amsfonts}
\usepackage{mathrsfs}

\usepackage[mathcal]{euscript}

\usepackage[ulem=normalem,draft]{changes}




%
\newtheorem{theorem}{Theorem}[section]

\newtheorem{corollary}[theorem]{Corollary}

\newtheorem{lemma}[theorem]{Lemma}

\finqui

\def\Beq{\Begin{equation}}
\def\Eeq{\End{equation}}

\def\Bthm{\Begin{theorem}}
\def\Ethm{\End{theorem}}
\def\Blem{\Begin{lemma}}
\def\Elem{\End{lemma}}

\def\Brem{\Begin{remark}\rm}
\def\Erem{\End{remark}}

\def\Bdim{\Begin{proof}}
\def\Edim{\End{proof}}




\def\step #1 \par{\medskip\noindent{\bf #1.}\quad}


\def\lhs{left-hand side}
\def\rhs{right-hand side}



\def\multibold #1{\def\arg{#1}%
  \ifx\arg\pto \let\next\relax
  \else
  \def\next{\expandafter
    \def\csname #1#1#1\endcsname{{\bf #1}}%
    \multibold}%
  \fi \next}

\def\pto{.}

\def\multical #1{\def\arg{#1}%
  \ifx\arg\pto \let\next\relax
  \else
  \def\next{\expandafter
    \def\csname cal#1\endcsname{{\cal #1}}%
    \multical}%
  \fi \next}


\def\multimathop #1 {\def\arg{#1}%
  \ifx\arg\pto \let\next\relax
  \else
  \def\next{\expandafter
    \def\csname #1\endcsname{\mathop{\rm #1}\nolimits}%
    \multimathop}%
  \fi \next}

\multibold
qwertyuiopasdfghjklzxcvbnmQWERTYUIOPASDFGHJKLZXCVBNM.

\multical
QWERTYUIOPASDFGHJKLZXCVBNM.

\multimathop
dist div dom meas sign supp .


\def\Accorpa #1#2 #3 {\gdef #1{\eqref{#2}--\eqref{#3}}%
  \wlog{}\wlog{\string #1 -> #2 - #3}\wlog{}}


\def\<#1>{\mathopen\langle #1\mathclose\rangle}
\def\norma #1{\mathopen \| #1\mathclose \|}


\def\iO{\int_\Omega}
\def\iQt{\iint_{Q_t}}
\def\iQ{\iint_Q}

\def\dt{\partial_t}
\def\dn{\partial_{\bf n}}

\def\checkmmode #1{\relax\ifmmode\hbox{#1}\else{#1}\fi}


\def\erre{{\mathbb{R}}}
\def\enne{{\mathbb{N}}}

\def\oT{(0,T)}





\def\Lx #1{L^{#1}(\Omega)}
\def\Hx #1{H^{#1}(\Omega)}

\def\Ldue{\Lx 2}

\def\Huno{\Hx 1}
\def\Hdue{\Hx 2}

\def\LiQ{L^\infty(Q)}
\def\LiS{L^\infty(\Sigma)}



\let\theta\vartheta

\let\phi\varphi

\let\TeXchi\chi                         
\newbox\chibox
\setbox0 \hbox{\mathsurround0pt $\TeXchi$}
\setbox\chibox \hbox{\raise\dp0 \box 0 }
\def\chi{\copy\chibox}


\def\pG{p_\Gamma}
\def\aG{\alpha_\Gamma}
\def\vG{v_\Gamma}
\def\wG{w_\Gamma}
\def\fG{f_\Gamma}
\def\feG{f_{\Gamma,1}}
\def\fzG{f_{\Gamma,2}}

\def\uG{u_\Gamma}
\def\yG{y_\Gamma}

\def\nabG{\nabla_\Gamma}
\def\delG{\Delta_\Gamma}
\def\iG{\int_\Gamma}
\def\iS{\iint_\Sigma}
\def\iSt{\iint_{\Sigma_t}}

\def\HG{{H_\Gamma}}
\def\VG{{V_\Gamma}}
\def\CV{{\cal V}}
\def\CH{{\cal H}}
\def\CS{{\cal S}}
\def\CP{{\bf {(CP)}}}

\def\CY{{\cal Y}}

\def\CU{{\cal U}}
\def\WG{{W_\Gamma}}
\def\CW{{\cal W}}
\def\Uad{{\cal U}_{\rm ad}}
\def\CZ{{\cal Z}}
\def\CN{{\cal N}}

\def\rhomin{\rho_{\rm min}}
\def\Rmin{\rho_{\Gamma_{\rm min}}}
\def\rhomax{\rho_{\rm max}}
\def\Rmax{\rho_{\Gamma_{\rm max}}}
\def\pG{p_\Gamma}
\def\bl{{\boldsymbol \lambda}}

\def\us{u^*}
\def\ys{y^*}
\def\ysG{y_\Gamma^*}
\def\usG{u_\Gamma^*}
\def\ps{p^*}
\def\psG{p_\Gamma^*}
\def\xiG{\xi_\Gamma}
\def\etaG{\eta_\Gamma}
\def\hG{h_\Gamma}
\def\kG{k_\Gamma}

\def\bu{{\bf u}}
\def\bus{{\bf u^*}}
\def\by{{\bf y}}
\def\bys{{\bf y^*}}
\def\bv{{\bf v}}
\def\bal{{\boldsymbol \alpha}}
\def\bp{{\bf p}}
\def\bps{{\bf p^*}}
\def\bz{{\bf z}}

\def\wJ{{\widehat{\cal J}}}

\def\VD{V^*}


\let\tau \tau
\normalfont\let\hat\widehat

\Begin{document}


\title{{\bf Second-order sufficient conditions for sparse optimal control of singular Allen--Cahn systems with dynamic boundary conditions}}
\author{}
\date{}
\maketitle

\Bcenter
{\large\bf J\"urgen Sprekels$^{(1)}$}\\
{\normalsize e-mail: {\tt juergen.sprekels@wias-berlin.de}}\\[.4cm]
{\large\bf Fredi Tr\"oltzsch$^{(2)}$}\\
{\normalsize e-mail: {\tt  troeltzsch@math.tu-berlin.de}}\\[.6cm]
$^{(1)}$
{\small Weierstrass Institute for Applied Analysis and Stochastics}\\
{\small Mohrenstra\ss e 39, 10117 Berlin, Germany}\\[.2cm]
$^{(2)}$
{\small Institut f\"ur Mathematik der Technischen Universit\"at Berlin}\\
{\small Stra\ss e des 17. Juni 136, 10623 Berlin, Germany}\\[.8cm]
\Ecenter

{
\Begin{abstract}\noindent
In this paper we study the optimal control of a parabolic initial-boundary value problem of Allen--Cahn type with dynamic boundary
conditions. Phase field systems of this type govern the evolution of
{coupled diffusive phase transition processes with nonconserved order parameters 
that occur in a container and on its surface, respectively. } It is assumed that the nonlinear functions driving the physical processes within the bulk and on the surface are double well potentials of logarithmic type whose derivatives become singular at the boundary of their respective domains of definition. 
For such systems, optimal control problems have been studied in the past.  We focus here on the situation when the cost functional 
of the optimal control problem contains a nondifferentiable term like the $L^1$-norm leading to sparsity of optimal controls. 
For such cases, we derive second-order sufficient conditions for locally optimal controls.  
\\[2mm]
{\bf Key words:}
Allen--Cahn equation, phase field model, dynamic boundary condition, singular potential, optimal control, sparsity, optimality conditions.
\normalfont
\\[2mm]
\noindent {\bf AMS (MOS) Subject Classification:} 35K20, 35K55, 49J50,  49J52, 49K20.
\End{abstract}
}
\salta

\pagestyle{myheadings}
\newcommand\testopari{\sc Sprekels \ --- \ Tr\"oltzsch}
\newcommand\testodispari{\sc Second-order conditions with sparsity for an Allen--Cahn system}
\markboth{\testodispari}{\testopari} 


\finqui


\section{Introduction}
\label{Intro}
\setcounter{equation}{0}

Let $\Omega\subset \erre^3$ denote some bounded and connected open set with smooth boundary $ \Gamma=\partial\Omega$ (a 
compact hypersurface of class $C^2$) and unit outward 
normal ${\bf n}$, and let, in this order, $\dn$, $\nabG$, $\delG$ denote the outward normal derivative, the surface gradient, and the 
Laplace--Beltrami operator 
defined on $\Gamma$. 
Recall that for sufficiently smooth functions $v$ defined on $\overline\Omega$ we have, with $\vG:=v_{|\Gamma}$,  the identities
$$\nabG \vG=\nabla v-(\nabla v\cdot {\bf n}){\bf n}, \quad \delG\vG=\nabG\cdot\nabG \vG,\quad\mbox{on }\,\Gamma.$$
Moreover, let $T>0$ denote a fixed final time, and 
\begin{align*}
&Q_t:=\Omega\times (0,t),\quad \Sigma_t:=\Gamma\times (0,t), \quad \mbox{for }\,t\in(0,T],\quad\mbox{and}\quad
Q:=Q_T,\quad \Sigma:=\Sigma_T.
\end{align*}
We then study optimal control problems for the parabolic system with nonlinear dynamic boundary condition
\begin{align}
\label{ss1}
&\dt y-\Delta y+f'(y)=u &&\quad\mbox{a.e. in }\,Q,\\
\label{ss2}
&\dt \yG-\delG \yG+\dn y+\fG'(\yG)=\uG \quad\mbox{and}\quad \yG=y_{|\Gamma} &&\quad\mbox{a.e. on }\,\Sigma,\\
\label{ss3}
&y(0)=y_0\,\mbox{ in }\,\Omega,\quad\yG(0)=y_{0_\Gamma}\quad\mbox{on }\,\Gamma.&&
\end{align}
\Accorpa\State ss1 ss3
Here, the second identity in \eqref{ss2} has to be interpreted in the following way: for almost every $t\in\oT$ the trace $y(t)_{|\Gamma}$ of
$y(t)$ on the boundary $\Gamma$ coincides with $y_\Gamma(t)$. Moreover, $y_0\,$ and $\, y_{0_\Gamma}\,$ are given initial data
with $y_{0_{|\Gamma}}=y_{0_\Gamma}$, and $\,u\,$ and $\, \uG\,$ are distributed and boundary controls, respectively. The pair $(y,y_\Gamma)$
is the state associated with the control pair $(u,u_\Gamma)$.
 
The system \State\ constitutes a phase field model for the physical process when an isothermal diffusive phase transition with 
nonconserved order parameter $\,y\,$  taking place in the interior of a container $\Omega\subset\erre^3$ is coupled
via the dynamic boundary condition \eqref{ss2} to another nonconserving diffusive phase transition occurring on the 
surface $\Gamma$ of the container. We assume here that the order parameter $\,y\,$ (and thus also $\yG$) is normalized 
to attain its values 
in the interval $[-1,1]$, where the level sets $\{y=-1\}$ and $\{y=1\}$ correspond to the pure phases. 
We remark that in our setting the pure phases will never occur (see the separation property \eqref{separation} below).
For more information on
the physical background of this model, we refer to \cite{SW} and the literature cited therein. 

A very important role for the evolution play the nonlinear functions $\,f\,$ and $\,\fG$, which are double well potentials whose
derivatives define the thermodynamic forces driving the phase transitions in the bulk and on the surface, respectively. Typical physically
relevant cases are given by the {\em regular}, {\em logarithmic}, and {\em double obstacle} potentials. In this order,  they 
 are given by
\begin{align}
\label{freg}
f_{\rm reg}(r)&=\frac 14(r-1)^2 \quad\mbox{for }\,r\in\erre,\\
\label{flog}
f_{\rm log}(r)&=\left\{\begin{array}{ll}
c_1\bigl((1+r)\,\ln(1+r)+(1-r)\,\ln(1-r)\bigr)-c_2r^2 &\quad\mbox{if }\,r\in(-1,1)\\
2c_1\,\ln(2)-c_2 &\quad\mbox{if }\,r\in\{-1,1\}\\
+\infty&\quad\mbox{if }\,r\not\in [-1,1]
\end{array}
\right. \,,\\
\label{f2obs}
f_{\rm 2obs}(r)&=\left\{\begin{array}{ll}
-c_3r^2&\quad\mbox{if }\,r\in[-1,1]\\
+\infty&\quad\mbox{if }\,r\not\in[-1,1]
\end{array}
\right.\,. 
\end{align}
In this connection, we have $c_i>0$, $i=1,2,3$, and $c_1,c_2$ are such that $f_{\rm log}$ is nonconvex.

There exists a vast literature on the well-posedness and asymptotic behavior of the Allen--Cahn equation and its various generalizations 
when complemented 
with dynamic boundary conditions. Without claiming to be comprehensive, we refer the reader to the works 
\cite{CC,ChGrM,Chill,CoFu, CGNS,GGM1,GGM2,Israel,Liero,SW,Wu}. Optimal control problems for \State\ have been studied in \cite{BEM} for the regular case and
in \cite{CS} for the singular logarithmic case \eqref{flog}, while the double obstacle case \eqref{f2obs} was investigated in \cite{CFS}. In particular, first-order necessary
and second-order sufficient optimality conditions have been derived in \cite{CS} for the differentiable tracking-type cost functional
\begin{align}
\label{diffpart}
J((y,\yG),(u,\uG))\,:=\,&\frac{\beta_1}2\iQ|y-y_Q|^2\,+\,\frac{\beta_2}2\iS|\yG-y_\Sigma|^2 \,+\,\frac{\beta_3}2
\iO|y(T)-y_{\Omega,T}|^2\nonumber\\
&+\,\frac{\beta_4}2\iG|\yG(T)-y_{\Gamma,T}|^2\,+\,\frac{\nu}2\iQ|u|^2\,+\,\frac{\nu_\Gamma }2\iS|\uG|^2\,,
\end{align}
which was also considered in \cite{CFS}. Here, the targets $y_Q, y_\Sigma, y_{\Omega,T},y_{\Gamma,T}$ are given functions,
$\beta_i$, $1\le i\le 4$, denote nonnegative constants which are not simultaneously zero, and \,$\nu,\nu_\Gamma\,$ are
positive constants. 

In this paper, we focus on the aspect of sparsity. To this end, we add to the cost functional $\,J\,$ a suitable term forcing sparsity, 
which is typically of the form
\begin{align} 
\label{nondiff}
j((u,\uG))\,:=\,\alpha\iQ|u|\,+\,\alpha_\Gamma\iS|\uG|,
\end{align}
with positive coefficients $\alpha,\aG$. The total cost functional is then given by 
\begin{equation}
\label{cost}
{\cal J}((y,\yG),(u,\uG))\,:=\,J((y,\yG),(u,\uG))+j((u,\uG))\,.
\end{equation}
Notice that $\,j\,$ is nonsmooth, which then also applies to ${\cal J}$. 

At this point, we specify the set of admissible controls: choosing real constants $\rhomin$, $\rhomax$, $\Rmin$, $\Rmax$ with $\rhomin\le\rhomax$ and
$\Rmin\le\Rmax$, we set
\begin{align}
\label{Uad}
\Uad\,:=&\,\left\{(u,\uG)\in L^\infty(Q)\times L^\infty(\Sigma):\rhomin\le u\le\rhomax
\,\mbox{ \,\,a.e. in }\,Q, \right.
\nonumber\\
&\left. \quad
\Rmin\le\uG\le\Rmax \mbox{ \,\,a.e. on }\,\Sigma\right\}.
\end{align}
\Brem
We remark that the above constants could also be replaced by functions belonging to $L^\infty(Q)$ and $L^\infty(\Sigma)$, respectively. 
We also stress the fact that, in contrast to the state functions $y$ and $y_\Gamma$,
the controls $u$ and $\uG$ are completely independent from each other and not linked on the
boundary. 
We can therefore treat here the cases 
of distributed and of boundary controls simultaneously; the cases of controling either only in the bulk or only on the surface are  
obtained by putting some of the weights $\beta_i$ and the associated thresholds defining $\Uad$ equal to zero correspondingly.  
\Erem

The control problem under investigation in this paper now reads as follows:

\vspace{2mm}\noindent
\CP \quad Minimize $\,\,{\cal J}((y,\yG),(u,\uG))\,\,$ subject to \State\ and $(u,\uG)\in \Uad.$

\vspace{2mm}
There are numerous papers studying optimal control problems for problems with dynamic boundary conditions. Without claiming to be
complete, we cite here the works \cite{BEM,CFS,CS} for the Allen--Cahn equation and 
\cite{CFGS1,CFGS2,CGSANA,CGSAMO, CGSAnnali, CGSSIAM, CGSconv,CSig,GS}
for systems of Cahn--Hilliard type. 

Sparsity in the optimal control theory for partial differential equations has become a very active field of research.  
The use of sparsity-enhancing functionals goes back to inverse problems and image processing. It was the seminal paper  
\cite{stadler2009}  on elliptic control problems that initiated the discussion of sparsity in the optimal control theory 
of partial differential equations. Soon after  \cite{stadler2009}, many results on sparse optimal controls for PDEs 
were published. We mention only very few of them with closer relation to our paper, in particular 
\cite{casas_herzog_wachsmuth2017,herzog_obermeier_wachsmuth2015,herzog_stadler_wachsmuth2012}, on directional sparsity,
 and \cite{casas_troeltzsch2012} on a general theorem for second-order conditions. 
Moreover, we refer to some new trends in the investigation of sparsity, namely, infinite horizon sparse 
optimal control (see, e.g., \cite{Kalise_Kunisch_Rao2017,Kalise_Kunisch_Rao2020}) and fractional order optimal 
control  (cf. \cite{Otarola_Salgado2018}, \cite{Otarola2020}).

These papers concentrated on first-order optimality conditions for sparse optimal controls of single elliptic and parabolic equations. 
In  \cite{casas_ryll_troeltzsch2013,casas_ryll_troeltzsch2015}, first- and second-order optimality conditions 
have been discussed in the context of sparsity for the (semilinear) system of  FitzHugh--Nagumo equations. 
More recently, sparsity of optimal controls for reaction-diffusion systems of Cahn--Hilliard type have been addressed in 
\cite{CSS4,Garcke_etal2021,SpTr1}.  Moreover, we refer to the measure control of the Navier--Stokes system 
studied in \cite{Casas_Kunisch2021}.

However, to the best knowledge of the authors, second-order sufficient optimality for sparse controls 
for the Allen--Cahn equation with dynamic boundary conditions have never been studied before. We also improve 
a result on second-order sufficient conditions of \cite{CS} for the case without sparsity functionals.

The paper is organized as follows. 
In the next section, we list and discuss our assumptions, and we collect known results concerning the well-posedness of the state system \State\ and of the optimal control problem that 
have been established in \cite{CC} and \cite{CS}.

In Section 3, we employ the implicit function 
theorem to give a new proof of the known fact that the control-to-state operator $(u,\uG)\mapsto (y,\yG)$ is twice continuously Fr\'echet differentiable between appropriate Banach spaces. The final Section 4 then brings the main results of this paper, 
namely the derivation of first-order necessary and second-order sufficient optimality conditions for the optimal control problem \CP. In an appendix, we prove auxiliary results that are needed for the main theorem on second-order sufficient conditions.

Prior to this, let us fix some notation.
For any Banach space $X$, we denote by \,$\|\,\cdot\,\|_X$, $X^*$, and $\langle\, \cdot\, , \,\cdot\, \rangle_X$,  
the corresponding norm, its dual space, and  the related duality pairing between $X^*$ and~$X$. 
For two Banach spaces $X$ and $Y$ that are both continuously embedded in some topological vector space~$Z$, we introduce the linear space
$X\cap Y$ that becomes a Banach space when equipped with its natural norm 
$\norma v_{X\cap Y}:=\norma v_X+\norma v_Y\,$  for $v\in X\cap Y$. 

The standard Lebesgue and Sobolev spaces defined on a set $G$ (where here $G=\Omega$ or  $G=\Gamma$) are, for
$1\le p\le\infty$ and $k\ge 0$, denoted by $L^p(G)$ and $W^{k,p}(G)$, respectively.  
If $p=2$, they become Hilbert spaces, and we employ the standard convention $H^k(\Omega):= W^{k,2}(\Omega)$. 
For convenience, we also~introduce the notation
\begin{align*}
  & H := \Ldue , \quad \HG:=L^2(\Gamma) ,\quad \CH:=H\times\HG,\\
  &V := \Huno,   \quad\VG:=H^1(\Gamma),\quad\CV:=\{(v,\vG)\in V\times\VG: \vG=v_{|\Gamma}\},\\
	&W:=\Hdue, \quad \WG:=H^2(\Gamma), \quad \CW:=W\times\WG.
  \end{align*}
All of these spaces are Banach spaces when endowed with their natural norms.	
We denote by $(\,\cdot\,,\,\cdot\,)_H$ the natural inner product in $H$. 
As usual, $H$ is identified with a subspace of $\VD$ according to the identity
\begin{align*}
	\langle u,v\rangle_V =(u,v)_H
	\quad\mbox{for every $u\in H$ and $v\in V$}.
\end{align*}
We then have the Hilbert triple $(V,H,\VD)$ with dense and compact embeddings. In the same way, using the canonical inner products
$$
(\vG,\wG)_{\HG}=\int_\Gamma \vG \wG \quad\mbox{and}\quad
((v,\vG),(w,\wG))_{\CH}=\iO vw + \iG \vG \wG
$$
in the spaces $\HG$ and $\CH$, respectively, we can construct the Hilbert triples  
$(\VG,\HG,\VG^*)$ and $(\CV,\CH,\CV^*)$ with dense and compact embeddings. We also recall the following version of Green's
formula for functions $\vG\in H^2(\Gamma)$ and $\wG\in H^1(\Gamma)$, which is valid since the compact hypersurface $\Gamma$ has 
an empty boundary (for details, see, e.g., \cite{DE}):
\begin{equation}
\label{GreenG}
-\int_\Gamma \wG\,\delG\vG\,=\,\int_\Gamma \nabG\wG\cdot\nabG\vG\,.
\end{equation}

We close this section by introducing a convention concerning the constants used in estimates within this paper: we denote by $\,C\,$ any 
positive constant that depends only on the given data occurring in the state system and in the cost functional, as well as 
on a constant that bounds the $\left(L^\infty(Q)\times L^\infty(\Sigma)\right)$--norms of the elements of $\Uad$. The actual value of 
such generic constants $\,C\,$ 	may 
change from formula to formula or even within formulas. Finally, the notation $C_\delta$ indicates a positive constant that
additionally depends on the quantity $\delta$.   


\section{General assumptions and the state system}
\setcounter{equation}{0}

In this section, we formulate the general assumptions for the data of the state system \State, and we collect some known results for
the state system. Throughout this paper, we generally assume:
\begin{description}
\item[(A1)] \,\,$f=f_1+f_2$ and $\fG=\feG+\fzG$, where $f_1,\feG:\erre\to [0,+\infty]$ are lower semicontinuous and
convex with $f_1(0)=0$ and $\feG(0)=0$. Moreover, $f_2,\fzG:\erre\to\erre$ have Lipschitz continuous first derivatives on $\erre$.
In addition, we require  $f_2, \fzG\in C^4[-1,1]$ and \,$f_1,\feG\in C^4(-1,1)$, and assume that
\begin{align}
\label{limf1}
&\lim_{r\searrow -1} \,f_1'(r)=\lim_{r\searrow -1} \,\feG'(r)=-\infty\,,\quad
\lim_{r\nearrow 1} \,f_1'(r)=\lim_{r\nearrow 1} \,\feG'(r)=+\infty\,,\\[2mm]
\label{dominiert}
&\exists \,\,M_1\ge 0, \,M_2>0\,\mbox{ such that }\,\,\,|f_1'(r)|\le 
M_1+M_2\,|\feG'(r)|\,\quad\forall\,r\in (-1,1).\quad
\end{align}
\item[(A2)] \,\,$(y_0,y_{0_\Gamma})\in \CV \cap (L^\infty(\Omega)\times L^\infty(\Gamma))$, and it holds that
\begin{align}
\label{ini1}
&-1\,<\, {\rm ess\,inf}_{x\in\Omega}\, y_0(x), \quad \,\,\,{\rm ess\,sup}_{x\in\Omega}\,y_0(x)<1\,.
\end{align}
\item[(A3)] \,\,$R>0$  is a fixed constant such that 
\Beq
\label{defUR}
\Uad\subset   \CU_R :=\{(u,\uG)\in L^\infty(Q)\times L^\infty(\Sigma):
\,\,\|u\|_{L^\infty(Q)}\,+\,\|\uG\|_{L^\infty(\Sigma)}\,<R\}.
\Eeq  
\end{description} 

\vspace*{1mm}
\Brem
We observe that the condition {\bf (A1)} is fulfilled if both $\,f\,$ and $\,\fG\,$ are given by logarithmic expressions 
of the same type as  the potential $f_{\rm log}$ defined in \eqref{flog}. The condition \eqref{dominiert} is of
technical nature; it is needed in the proof of the existence result of Theorem 2.2 below. From the viewpoint of physics, it means that the thermodynamic force exerted on the surface somehow dominates the one acting in the bulk. Moreover, the condition \eqref{ini1} implies
that
$$
-1\,<\,{\rm ess\,inf}_{x\in\Gamma}\,y_{0_\Gamma}(x), \quad {\rm ess\,sup  }_{x\in\Gamma}\,y_{0_\Gamma}(x)\,<\,1.
$$
Therefore, the state $(y,\yG)$ is initially strictly separated from the endpoints of the interval $(-1,1)$. This means that
initially there are no pure phases within the container and on its surface. Finally, the condition {\bf (A3)} just fixes once and for all
a bounded open subset of the control space $\LiQ\times\LiS$ that contains $\Uad$. 
\Erem

Next, we specify our notion of solution: for any given $(u,\uG)\in \CH$, we call a pair $(y,\yG)$ a solution
to \State\  if 
\begin{align*}
& y \in  H^1(0,T;V^*)\cap L^2(0,T;V), \\
& \yG\in  H^1(0,T;V^*_{\Gamma})\cap L^2(0,T;\VG),\\
&\yG(t)=y(t)_{|\Gamma} \,\,\mbox{ for a.e. }\,t\in\oT,\\
& y(0)=y_0,\quad \yG(0)=y_{0_\Gamma},
\end{align*}
and if, for almost every $t\in (0,T)$ and every $(v,\vG)\in\mathcal{V}$, it holds
\begin{align}
\label{ssvar}
& \langle\partial_t y(t),v\rangle_V+\int_{\Omega}\nabla y(t)\cdot\nabla v + \langle\partial_t \yG(t),\vG\rangle_{V_{\Gamma}} +
\int_{\Gamma}\nabla_{\Gamma}\yG(t)\cdot\nabla_{\Gamma}\vG \nonumber\\
&=\int_{\Omega}(u(t)-f'(y(t)))v+\int_{\Gamma}(\uG(t)-f'_{\Gamma}(\yG(t)))\vG \,.
\end{align} 
Note that the identity \eqref{ssvar} is formally derived as follows: test \eqref{ss1} by $v$ and \eqref{ss2} by $\vG$, integrate by parts 
using \eqref{GreenG}, and add the resulting identities.

We have the following well-posedness result that follows from \cite[Thm.~2.1~and~Lem. 2.3]{CS}.
\Bthm
Suppose that the conditions {\bf (A1)}--{\bf (A3)} are fulfilled. Then the state system \State\ has for any 
$(u,\uG)\in L^2(Q)\times L^2(\Sigma)$ a unique solution $(y,\yG)$
with the regularity
\begin{align}
\label{reguy}
&y \in H^1(0,T;H)\cap L^\infty(0,T;V) \cap L^2(0,T;W),\\
\label{regyG}
&\yG\in H^1(0,T;\HG)\cap L^\infty(0,T;\VG)\cap L^2(0,T;\WG).
\end{align}
Moreover, there is a constant $K_1>0$, which depends only on $\,R\,$ and the data of the state system, such that 
\begin{align}
\label{ssbound1}
&\|y\|_{H^1(0,T;H)\cap L^\infty(0,T;V) \cap L^2(0,T;W)} \nonumber\\
&+\,\|\yG\|_{H^1(0,T;\HG)\cap L^\infty(0,T;\VG)\cap L^2(0,T;\WG)}\,\le\,K_1\,,
\end{align} 
whenever $(y,\yG)$ is the solution to the state system associated with some $(u,\uG)\in\CU_R$. In addition, a uniform strict separation
property is satisfied: there are constants $-1<r_- \le r+<1$, which depend only on $\,R\,$ and the data of the state system,
such that
\Beq
\label{separation}
r_- \le y(x,t) \le r_+ \,\,\mbox{ for a.e. }(x,t)\in Q,  \quad r_- \le \yG(x,t) \le r_+ \,\,\mbox{for a.e. }(x,t)\in\Sigma,
\Eeq
whenever $(y,\yG)$ is the solution to the state system associated with some $(u,\uG)\in\CU_R$.
\Ethm
\Brem
1. By virtue of Theorem 2.2, the control-to-state operator $\CS:(u,\uG)\mapsto \CS(u,\uG):=(y,\yG)$ is well defined as a mapping
between $L^2(Q) \times L^2(\Sigma)$ and the Banach space 
$H^1(0,T;\CH)\cap L^\infty(0,T;\CV)\cap L^2(0,T;\CW)$,
which also encodes 
the condition that
$\yG(t)=y(t)_{|\Gamma}$ for almost every $t\in\oT$. In particular, $(y,\yG)$ is a strong solution
to the state system that satisfies the equations \State\ almost everywhere.\\
2. Observe that the separation condition \eqref{separation} holds only for (bounded) controls in $\CU_R$. If it is satisfied, then,
by condition {\bf (A1)},  we may without loss of generality assume that
\Beq
\label{ssbound2}
\max_{i=1,2,3,4}\,\max_{j=1,2}\,\left(\|f_j^{(i)}(y)\|_{\LiQ}\,+\,\|f_{\Gamma,j}^{(i)}(\yG)\|_{\LiS}\right)\,\le\,K_1\,,
\Eeq
for every solution $(y,\yG)$ associated with some $(u, \uG)\in\CU_R$.\\
3. We cannot expect $\,y\,$ to be continuous on $\overline Q$, in general. 
However, we have that $y\in L^2(0,T;C^0(\overline\Omega))$ by the embedding $\Hdue\subset C^0(\overline\Omega)$.
This fact justifies our denotation for the trace on $\Gamma$: indeed, we have for almost all $t\in\oT$ that $y(t)\in V\cap 
C^0(\overline\Omega)$, and therefore the trace $y_\Gamma(t)$  coincides with the restriction of $y(t)$ to the boundary.\\
4. Since the embedding $\,(H^1(0,T;\CH)\cap L^2(0,T;\CV\cap\CW))\subset C^0([0,T];\CV)\,$ is continuous,  
the terminal observation $(y(T),y_\Gamma(T))$  in the functional \eqref{diffpart} is well defined.
\Erem


\section{Differentiability of the control-to-state operator}
\setcounter{equation}{0}

In this section, we study the differentiability properties of the control-to-state operator $\CS$. To this end, we introduce the
Banach spaces
\begin{align}
\label{defU}
&\CU:=\LiQ\times\LiS,\\
\label{defY}
&\CY:=H^1(0,T;\CH)\cap L^\infty(0,T;\CV)\cap L^2(0,T;\CW),
\end{align}
endowed with their standard norms.  
We then know from Theorem 3.2 and Theorem 3.5 in \cite{CS}
that $\CS$ is under the assumptions {\bf (A1)}--{\bf (A3)} twice continuously Fr\'echet differentiable on $\CU$ 
as a mapping from $\CU$ into $\CY$, where,  for any control pair $(\us,\usG)\in\CU$, with 
associated state $(\ys,\ysG):=\CS(\us,\usG)$, the
first and second Fr\'echet derivatives $D\CS(\us,\usG)\in {\cal L}(\CU,\CY)$ and $D^2\CS(\us,\usG)\in{\cal L}(\CU,{\cal L}(\CU,\CY))$
are given as follows:

\vspace{1mm}\noindent
(i) \,\,\,\,For any  increment $(h,\hG)\in\CU$, $(\xi,\xiG):=D\CS(\us,\usG)[(h,\hG)]\in\CY$ is the unique solution to the linearized problem
\begin{align}
\label{ls1}
&\dt\xi-\Delta\xi+f''(\ys)\xi=h\quad\mbox{a.e. in }\,Q,\\
\label{ls2}
&\dt\xiG-\delG\xiG+\dn\xi+f_{\Gamma}''(\ysG)\xiG=\hG\quad\mbox{and}\quad\xiG=\xi_{|\Gamma}\quad\mbox{a.e. on }\,\Sigma,\\
\label{ls3}
&\xi(0)=0 \,\mbox{\, a.e. in \,}\Omega,\quad \xiG(0)=0 \,\mbox{\, a.e. on \,}\Gamma.
\end{align} 
\Accorpa\Linear ls1 ls3

\vspace{1mm}
\noindent
(ii) \,\, For any pair of increments $(h,\hG),(k,\kG)\in \CU$, $(\eta,\etaG):=D^2\CS(\us,\usG)[(h,\hG),(k,\kG)]\in \CY$ is the unique solution to the
bilinearized problem
\begin{align}
\label{bilin1}
&\dt\eta-\Delta\eta+f''(\ys)\eta=-f^{(3)}(\ys)\phi\psi
\quad\mbox{a.e. in }\,Q,\\
\label{bilin2}
&\dt\etaG-\delG\etaG+\dn\eta + f_{\Gamma}''(\ysG)\etaG=-\fG^{(3)}(\ysG)\phi_\Gamma\psi_\Gamma
\quad\mbox{and}\quad \etaG=\eta_{|\Gamma}\quad\mbox{a.e. on  }\,\Sigma,\\
\label{bilin3}
&\eta(0)=0 \quad\mbox{a.e. in }\,\Omega, \quad \etaG(0)=0 \quad\mbox{a.e. on }\,\Gamma,
\end{align}
\Accorpa\Bilinear bilin1 bilin3
where $\,(\phi,\phi_\Gamma):=D\CS(\us,\usG)[(h,\hG)]\,$ and $\,(\psi,\psi_\Gamma):=D\CS(\us,\usG)[(k,\kG)]$.

\vspace{1mm}\noindent
(iii) \,The mappings $\,D\CS:\CU\to {\cal L}(\CU,\CY),(u,\uG) \mapsto D\CS(u,\uG),$\, and 
$\,D^2\CS(u,\uG):\CU\to {\cal L}(\CU,{\cal L}(\CU,\CY)), (u,\uG)\mapsto D^2\CS(u,\uG)$,\, are 
Lipschitz continuous in the following sense: there exists a constant
$K_2>0$, which depends only on $R$ and the data, such that, for all controls $(u,\uG), (\us,\usG)\in\CU_R$ and all 
increments $(h,\hG),(k,\kG)\in\CU$,
\begin{align}
\label{lip1} 
&\|(D\CS(u,\uG)-D\CS(\us,\usG))[(h,\hG)]\|_{\CY}\nonumber\\
&\quad\le\,K_2\,\|(u,\uG)-(\us,\usG)\|_{L^2(0,T;\CH)}\,\|(h,\hG)\|_{L^2(0,T;\CH)}\,,\\[2mm]
\label{lip2}
&\left\|\left(D^2 \CS  (u,\uG)-D^2 \CS (\us,\usG)\right)[(h,\hG),(k,\kG)]\right\|_{\CY}  \nonumber\\
&\quad\le\,K_2\,\|(u,\uG)-(\us,\usG)\|_{L^2(0,T;\CH)}\,\|(h,\hG)\|_{L^2(0,T;\CH)}\,\|(k,\kG)\|_{L^2(0,T;\CH)}\,.
\end{align}

\vspace*{3mm}
\Brem
As $\CU$ is dense in $L^2(0,T;\CH)$, the operator $D\CS(\us,\usG)\in {\cal L}(\CU,\CY)$ can be extended in the standard way to an operator 
belonging to ${\cal L}(L^2(0,T;\CH),\CY)$ without changing its operator norm. 
 We still denote the extended operator by $D\CS(\us,\usG)$, where we stress the fact that it coincides with a Fr\'echet derivative only on $\CU$ and not on $L^2(0,T;\CH)$. However, it follows from \cite[Thm.~2.2]{CS} that the linearized system \Linear\
has also for every \rhs\ $(h,\hG)\in L^2(0,T;\CH)$ a unique solution $(\xi,\xiG)\in\CY$ that satisfies
$$\|(\xi,\xiG)\|_{\CY}\,\le\,K_3\,\|(h,\hG)\|_{L^2(0,T;\CH)}$$
with a constant $K_3>0$ that depends only on $R$ and the data. It is then easily verified that $(\xi,\xiG)=D\CS(\us,\usG)[(h,\hG)]$ with the
extended operator, and, in the sense of the extension, the estimate \eqref{lip1} is also satisfied for directions $(h,\hG)\in L^2(0,T;\CH)$.  
An analogous result holds for the validity of \eqref{lip2}. 
\Erem

The above results (i)--(iii) have been proved directly in \cite{CS} without use of the implicit function theorem, where the authors 
announced that an alternative proof would be possible using the
implicit function theorem. This does not seem to be obvious, since the presence of  nonlinearities would require  
differentiability properties of Nemytskii operators between $L^\infty$--spaces. It is, however, not known whether the solutions to
linear systems like \Linear\ are bounded. Below (see Lemma 3.2), we will show such a boundedness result for bounded right-hand sides.
Using this result, we will be able to prove differentiability via the implicit function theorem.

To this end, we introduce the Banach space
\begin{align}
\label{defZ}
&\CZ:=\{(y,\yG)\in\CY\cap\CU: \,\dt y-\Delta y\in\LiQ,\,\,\dt\yG-\delG\yG+\dn y\in\LiS\},
\end{align}
endowed with the norm
\Beq
\label{normZ}
\|(y,\yG)\|_{\CZ}:=\|(y,\yG)\|_{\CY\cap\CU}+\|\dt y-\Delta y\|_{\LiQ}+\|\dt\yG-\delG\yG+\dn y\|_{\LiS}\,\,\,\forall (y,\yG)\in\CZ.
\Eeq
Finally, we fix constants $r_*,r^*$ such that
\Beq
\label{rstar}
-1<r_*<r_-<r_+<r^*<1,
\Eeq
with the constants $r_-,r_+$ introduced in \eqref{separation}. We then consider the set 
\begin{align}
\label{defPhi}
\Phi:=&\left\{(y,\yG)\in\CZ: r_* < \min\bigl\{{\rm ess\,inf}_{(x,t)\in Q}\,y(x,t)\,,\,{\rm ess\,inf}_{(x,t)\in\Sigma}\,\yG(x,t)\bigr\}\right.
\nonumber\\
&\qquad  \left. \mbox{and} \quad\max\bigl\{{\rm ess\,sup}_{(x,t)\in Q}\,y(x,t)\,,\,{\rm ess\,sup}_{(x,t)\in\Sigma}\,\yG(x,t)\bigr\}<r^* \right\},
\end{align}
which is obviously an open subset of $\CZ$. 
Notice that the functions in $\CZ$ are bounded and measurable, so that the essential infimum and supremum used above are well defined.

We now prove an auxiliary result for the linear initial-boundary value problem
\begin{align}
\label{aux1}
&\dt y-\Delta y\,=\,-\lambda_1\,f''(\ys)y+\lambda_2\,h \quad\mbox{a.e. in }\,Q,\\
\label{aux2}
&\dt \yG-\delG\yG+\dn y\,=\,-\lambda_1\,\fG''(\ysG)\yG + \lambda_2 \,\hG \quad\mbox{and} \quad \yG=y_{|\Gamma} \quad
\mbox{a.e. on }\,\Sigma,\\
\label{aux3}
&y(0)=\lambda_3 y_0 \quad\mbox{a.e. in }\,\Omega,\quad y_\Gamma(0)=\lambda_3 y_{0_\Gamma}\quad\mbox{a.e. on }\,\Gamma,
\end{align}
\Accorpa\Aux aux1  aux3
which for $\lambda_1=\lambda_2=1$ and $\lambda_3=0$ coincides with the linearization \Linear\ of the state system at
$((\us,\usG),(\ys,\ysG))$. For convenience, we now introduce the Banach space of the initial data,
\Beq
\label{defN}
{\cal N}:=\{(y_0,y_{0_\Gamma}):y_0\in V\cap L^\infty(\Omega), \,\,\,y_{0_\Gamma}\in \VG\cap L^\infty(\Gamma), \,\,\,
y_{0_\Gamma}=y_{0_{|\Gamma}}
\,\mbox{ a.e. on }\,\Gamma\},
\Eeq
equipped with its natural norm. We then have the following result.
\Blem
Assume that $\lambda_1,\lambda_2,\lambda_3 \in\{0,1\}$ are given and that the assumptions {\bf (A1)}--{\bf (A3)} are fulfilled.
Moreover, let $((\us,\usG),(\ys,\ysG))\in \CU_R\times\Phi$ be arbitrary. Then the system \Aux\ has for every $(h,\hG)\in\CU$ and
every $(y_0,y_{0_\Gamma})\in {\cal N}$ a unique solution $(y,\yG)\in\CZ$. Moreover, the linear mapping
$\,((h,\hG),(y_0,y_{0_\Gamma}))\mapsto (y,\yG)\,$
is continuous from $\,\CU\times{\cal N}\,$ into $\CZ$. 
\Elem
\Bdim
At first, it is standard to show that \Aux\ has a unique solution $(y,\yG)\in\CY$ for given data $(h,\hG)\in\CU$ and
$(y_0,y_{0_\Gamma})\in {\cal N}$. The existence can be proved via an appropriate Faedo--Galerkin 
approximation for which a priori estimates and a passage to the limit process are performed. The uniqueness proof is simple. 
In order not to overload the exposition, we 
avoid writing the Faedo--Galerkin scheme here and just give the corresponding a priori estimates formally.  To this end, we introduce the
constant
\Beq
\label{defM}
M:= \lambda_2\,\|(h,\hG)\|_{\CU}\,+\,\lambda_3\,\|(y_0,y_{0_\Gamma}\|_{\CN}.
\Eeq
Now we put
$$z:=-\lambda_1 f''(\ys)y+\lambda_2h,\quad z_\Gamma:=-\lambda_1 \fG''(\ysG)\yG+\lambda_2\hG.$$
Recalling that $(\ys,\ysG)\in \Phi$, and putting $\,\gamma:=1+\|f''\|_{C^0([r_*,r^*])}+\|\fG''\|_{C^0([r_*,r^*])}$, we have the estimates
\begin{align}
\label{rhs}
|z|\,\le\,\gamma(|y|+M) \quad\mbox{a.e. in }\,Q,\quad |z_\Gamma|\,\le\,\gamma(|\yG|+M)\quad\mbox{a.e. on }\,\Sigma .
\end{align}
In the remainder of the proof, we denote by $C>0$ constants that may depend on $\,\gamma\,$ but not on $\,M$.

Next, we add $y$ to both sides of \eqref{aux1} and $\yG$ to both sides of \eqref{aux2}, and we multiply the resulting identities by
$\dt y\in L^2(Q)$ and $\dt \yG\in L^2(\Sigma)$, respectively. Then we  integrate the results for arbitrary $t\in(0,T]$ over $Q_t$ and $\Sigma_t$, respectively, integrate by parts using \eqref{GreenG}, and add the results. We then arrive at the identity 
\begin{align}
\label{uwe1}
&\iQt|\dt y|^2 + \iSt|\dt\yG|^2 + \frac 12\|y(t)\|_V^2  + \frac 12 \|\yG(t)\|_{\VG}^2 \nonumber\\
&= \frac {\lambda_3^2}2 \|y_0\|_V^2 + \frac {\lambda_3^2} 2 \|y_{0_\Gamma}\|_{\VG}^2 +\iQt(y+z)\dt y + \iSt (\yG+z_\Gamma)\dt\yG\,.
\end{align}
Applying Young's inequality appropriately to the last two summands on the \rhs, using \eqref{rhs}, and then invoking Gronwall's lemma, we 
 easily conclude that
\begin{align}
\label{uwe2}
\|(y,\yG)\|_{H^1(0,T;\CH)\cap L^\infty(0,T;\CV)}\,\le\,C\,M.
\end{align} 

At this point, we observe that $(g,g_\Gamma):=(z-\dt y,z_\Gamma-\dt\yG)\in \CH$ almost everywhere in $\oT$. 
It therefore follows from the regularity result established in \cite[Lem.~3.1]{CGSAnnali} that for almost every $t\in\oT$
it holds $(y(t),\yG(t))\in\CW$, and, with a constant $C_\Omega>0$ that depends only on $\Omega$,
$$
\|(y(t),\yG(t))\|_{\CW}\,\le\,C_\Omega\left(\|(y(t),\yG(t))\|_{\CV}\,+\,\|(g(t),g_\Gamma(t))\|_{\CH}\right)\,. 
$$
Thus, using \eqref{uwe2} and \eqref{rhs}, we readily conclude that $(y,\yG)\in L^2(0,T;\CW)$ with
\Beq
\label{uwe3}
\|(y,\yG)\|_{L^2(0,T;\CW)}\,\le\,C\,M.
\Eeq
Combining \eqref{uwe2} and \eqref{uwe3}, we therefore have $(y,\yG)\in \CY$ with the bound
\Beq
\label{uwe4}
\|(y,\yG)\|_{\CY}\,\le\,C\,M.
\Eeq

Next, we are going to show that $(y,\yG)\in\CU$ with a corresponding norm estimate
\Beq
\label{uwe5}
\|(y,\yG\|_{\CU}\,\le\,C\,M.
\Eeq
 Once this will be shown,
\eqref{aux1}-\eqref{aux2} will yield that $(\dt y-\Delta y, \dt\yG-\delG\yG+\dn y)=(z,z_\Gamma) \in\CU$, which then implies that
$(y,\yG)\in \CZ$ with $\|(y,\yG)\|_{\CZ}\le C\,M$.

We argue by a Moser iteration technique. To this end, we rewrite the system \Aux. With the constant $\,\gamma\,$ introduced above, we 
put $w(x,t):=\exp (-\gamma t)y(x,t)$ and $w_\Gamma(x,t):=\exp(-\gamma t)\yG(x,t)$, noticing that for almost every $t\in\oT$ we have
$w_\Gamma(t)=w(t)_{|\Gamma}$. In terms of these new variables, the system \Aux\ becomes
\begin{align}
\label{hilf1}
&\dt w-\Delta w +(\gamma+f''(\ys))w= \exp(-\gamma t)\lambda_2h\quad\mbox{a.e. in }\,Q,
\\
\label{hilf2}
&\dt\wG-\delG\wG+\dn w+(\gamma+\fG''(\ysG))\wG= \exp(-\gamma t)\lambda_2\hG\quad\mbox{a.e. on }\,\Sigma,
\\
\label{hilf3}
&w(0)=\lambda_3 y_0\quad\mbox{a.e. in }\,\Omega, \quad w_\Gamma(0)=\lambda_3 y_{0_\Gamma}\quad\mbox{a.e. on }\,\Gamma.
\end{align}
\Accorpa\Hilfe hilf1 hilf3
We aim at showing that $\,(w,\wG)\in\CU\,$ and that, with a constant $\hat C>0$ not depending on $\,M$,
\Beq
\label{boundw}
\|(w,\wG)\|_{\CU}\,\le\,\hat C\,M.
\Eeq
Once this will be shown, we will have $\|(y,\yG)\|_{\CU}\,\le\,\exp(\gamma T)\,\hat C\,M$, and the proof of the assertion will be complete.
Observe that the system \Hilfe\ is suited better for proving an $L^\infty$--bound than \Aux, since the coefficient 
functions $c_0:=\gamma+f''(\ys)$ and $c_{0_\Gamma}:=\gamma+\fG''(\ysG)$ are nonnegative almost everywhere. In addition, the
\rhs s of \eqref{hilf1} and of \eqref{hilf2} are both bounded by $\,M$.

We now consider for $s>0$ the cutoff-functions
\Beq
\label{ws}
w^s:=\max\,\{-s,\min \{w,s\} \},\quad w_\Gamma^s:=\max\,\{-s,\min \{\wG,s\}\}\,.
\Eeq
We notice that $y\in L^2(0,T;C^0(\overline\Omega))$, by the embedding $\Hdue\subset C^0(\overline\Omega)$. Hence, for almost all
$t\in\oT$, we have $w^s(t)\in V\cap C^0(\overline\Omega)$, which means that the trace of $w^s(t)$ on $\Gamma$ is given by the
restriction of $w^s(t)$ to $\Gamma$.  In other words, it holds $w^s(t)_{|\Gamma}=w_\Gamma^s$, and $(w^s,w^s_\Gamma)\in \CV$. 
Moreover, we obviously have that $(w^s,w^s_\Gamma)\in\CU$. We 
therefore may for arbitrary integer $n\ge 2$ test the equations \eqref{hilf1} and \eqref{hilf2} by the admissible functions
$\,\,v=w^s(t)^{2n-1}\,\,$ and $\,\,\vG=w_\Gamma^s(t)^{2n-1}$, respectively. Integration by parts and over $[0,t]$, where $t\in (0,T]$,
addition of the resulting equalities, and the fact that $|e^{-\gamma t}|\le 1$, yield the inequality
\begin{align}
\label{uwe6}
&\iint_{Q_t}\dt w\,(w^s)^{2n-1}\,+\iint_{\Sigma_t}\dt\wG\,(w^s_\Gamma)^{2n-1} \,+\,(2n-1)\iQt |w^s|^{2n-2}\,
|\nabla w^s|^2 \nonumber\\
&+\,(2n-1)\iSt |w^s_\Gamma|^{2n-2}\,|\nabla_\Gamma\wG^s|^2\,+\,\iQt c_0 w (w^s)^{2n-1}\,+\,\iSt c_{0_\Gamma}\wG
(w^s_\Gamma)^{2n-1}\nonumber\\
&\le\,\iQt |\lambda_2h|\,|w^s|^{2n-1}  \, + \iSt |\lambda_2\hG|\,|w^s_\Gamma|^{2n-1}\,.
\end{align}

Now note that $2n-1$ is an odd integer, and thus the signs of $\,w$, $\,w^s\,$ and $\,(w^s)^{2n-1}\,$ are equal. But then, owing to the fact
that $c_0\ge 0$ by construction, the product $\,c_0 w(w^s)^{2n-1}\,$ is nonnegative almost everywhere. Hence, the fifth summand on the \lhs\
of \eqref{uwe6} is nonnegative, and, by the same token, also the sixth summand. 
Moreover, we have
\begin{align*}
&\iQt \dt w\,(w^s)^{2n-1}\,=\,\iQt \dt w^s\,(w^s)^{2n-1} \,+\iQt \dt(w-w^s)\,(w^s)^{2n-1} \nonumber\\
&=\,\frac 1{2n}\,\|w^s(t)\|_{L^{2n}(\Omega)}^{2n}-\frac 1{2n}\,\|w^s(0)\|_{L^{2n}(\Omega)}^{2n} \,+\iO (w(t)-w^s(t))(w^s(t))^{2n-1}\nonumber\\
&\quad - \iO(w(0)-w^s(0))(w(0))^{2n-1} \,-\iQt (w-w^s)\,\dt \bigl[(w^s)^{2n-1}\bigr] \,.
\end{align*}
Obviously, the integrand of the last summand on the \rhs\ is zero almost everywhere, and if we choose 
\begin{equation}
\label{erich}
s>\lambda_3\bigl(\|y_0\|_{L^\infty(\Omega)}+\|y_{0_\Gamma}\| _{L^\infty(\Gamma)}\bigr), 
\end{equation}
which will henceforth be assumed, then $w(0)=w^s(0)$, and also the integrand of the fourth summand on the \rhs\ vanishes. Finally,
we easily check that the integrand of the third term on the \rhs\ is nonnegative.
In summary, 
\begin{equation}
\label{uwe7}
\iint_{Q_t} \dt w\,(w^s)^{2n-1} \,\ge\, \frac 1{2n}\,\|w^s(t)\|_{L^{2n}(\Omega)}^{2n}-\frac 1{2n}\,\|w(0)\|_{L^{2n}(\Omega)}^{2n}
\,,                          
\end{equation}
and, by the same token, an analogous estimate holds true for the second summand on the left-hand side of \eqref{uwe6}. Hence, omitting several
nonnegative terms on the left-hand side of \eqref{uwe6},  we obtain from \eqref{uwe6}--\eqref{uwe7} the inequality
\begin{align}
\label{uwe8}
&\frac 1{2n}\,\|w^s(t)\|_{L^{2n}(\Omega)}^{2n}+\frac 1{2n}\,\|w^s_\Gamma(t)\|_{L^{2n}(\Gamma)}^{2n} \,\le\,
\frac 1{2n}\,\|w(0)\|_{L^{2n}(\Omega)}^{2n}+\frac 1{2n}\,\|w_\Gamma(0)\|_{L^{2n}(\Gamma)}^{2n} 
\nonumber\\
&\,+\iQt |\lambda_2 h|\,|w^s|^{2n-1}\,+\iSt |\lambda_2 \hG|\,|w^s_\Gamma|^{2n-1}\,.
\end{align}

It remains to estimate the terms on the \rhs. At first, denoting by $\kappa$ the maximum between the volume of $\Omega$ and 
the surface area of $\Gamma$, we have that
\begin{align}
\label{uwe9}                                                                                         
&\|w(0)\|_{L^{2n}(\Omega)}\,= \,\lambda_3\,\|y_0\|_{L^{2n}(\Omega)}\,\le\, \kappa^{1/(2n)}\,\lambda_3\,\|y_0\|_{L^\infty(\Omega)}
\,\le\, \kappa^{1/(2n)} \,M,\nonumber\\
&\|w_\Gamma(0)\|_{L^{2n}(\Gamma)}\,=\,\lambda_3\,\|y_{0_\Gamma}\|_{L^{2n}(\Gamma)}\,\le\,\kappa^{1/(2n)}\,
\lambda_3\,\|y_{0_\Gamma}\|_{L^\infty(\Gamma)}
\,\le\,\kappa^{1/(2n)}\,M.
\end{align}
In addition, we obtain for the third term on the \rhs\ (which we denote by $I$), using Young's inequality $\,\,ab\le \frac{1}p|a|^p+\frac 1q|b|^q\,\,$
with $\,p=\frac{2n}{2n-1}\,$ and $\,q=2n\,$,
\begin{align}
\label{uwe10}
I\,&\le \iQt M|w^s|^{2n-1}\,\le\,\frac{2n-1}{2n} \iQt|w^s|^{2n}\,+\,\frac{M^{2n}}{2n}\,|\Omega|\,t\nonumber\\
&\le \frac{2n-1}{2n} \iQt|w^s|^{2n}\,+\,\frac{M^{2n}}{2n}\,\kappa\,T.
\end{align}  
An analogous estimate can be performed for the last summand on the \rhs. Consequently, combining the estimates \eqref{uwe8}--\eqref{uwe10}, 
and multiplying the resulting inequality by $2n$, we arrive at the estimate
\begin{align*}
&\|w^s(t)\|_{L^{2n}(\Omega)}^{2n}\,+\, \|w^s_\Gamma(t)\|^{2n}_{L^{2n}(\Gamma)}\,\le\,
2\kappa(1+T)M^{2n}\\
&+\,(2n-1)\int_0^t\|w^s(\sigma)\|_{L^{2n}(\Omega)}^{2n}\,d\sigma\,+\,(2n-1)\int_0^t\|w^s_\Gamma(\sigma)\|^{2n}_{L^{2n}(\Gamma)}\,d\sigma\,,
\end{align*}
whence, by virtue of Gronwall's lemma,
\begin{align*}
&\|w^s(t)\|_{L^{2n}(\Omega)}^{2n}+ \|w^s_\Gamma(t)\|^{2n}_{L^{2n}(\Gamma)}\,\le\,
2\kappa(1+T)M^{2n}\,
e^{(2n-1)t}.
\end{align*}
Therefore,
\begin{align*}
&\|w^s(t)\|_{L^{2n}(\Omega)}\,\le \,\left(2\kappa(1+T)\right)^{1/(2n)}\,\exp(T)\,M,
\end{align*}
and, by the same token,
$$\|w^s_\Gamma(t)\|_{L^{2n}(\Gamma)}\,\le\,\left(2\kappa(1+T)\right)^{1/(2n)}\,\exp(T)\,M.$$
Taking the limit as $n\to\infty$ in the last two inequalities, we find that
\Beq
\label{uwe11} 
\|w^s(t)\|_{L^\infty(\Omega)} \,+\,\|w^s_\Gamma(t)\|_{L^\infty(\Gamma)}\,\le\,2\,\exp(T)\,M,
\Eeq
for almost every $t\in\oT$,  provided that $\,s\,$ satisfies \eqref{erich}.
At this point, we choose
$$s> 2\exp(T)M\,+\,\lambda_3\,\bigl(\|y_0\|_{L^\infty(\Omega)}\,+\,\|y_{0_\Gamma}\|_{L^\infty(\Gamma)}\bigr).$$
Then, owing to \eqref{uwe11}, $\,w^s=w\,$ and $\,w^s_\Gamma=\wG$, whence we conclude
that, for almost every $t\in\oT$, 
$$\|w(t)\|_{L^\infty(\Omega)}\,+\,\|\wG(t)\|_{L^\infty(\Gamma)}\,\le\,2\,\exp(T)\,M.$$
Hence \eqref{boundw} is shown, which concludes the proof of the assertion.   
\Edim

Having proved Lemma~3.2, we can now prepare for the application of the implicit function
theorem. To this end, we introduce for convenience abbreviating denotations, namely,
\begin{align*}
{\bf u}&:=(u,\uG),\quad {\bf u^*}:=(\us,\usG), \quad {\bf y}:=(y,\yG),\quad
{\bf y^*}:=(\ys,\ysG),\nonumber\\
{\bf y}_0&:=(y_0,y_{0_\Gamma}), \quad \mathbf{0}:=(0,0).
\end{align*} 
We consider two auxiliary linear initial-boundary value problems. The first,
\begin{align}
\label{sysG11} 
&\dt y-\Delta y\,= \, h\quad \mbox{a.e. in \,$Q$},\\
\label{sysG12} 
&\dt\yG-\Delta_\Gamma\yG+\dn y\,= \, \hG \quad\mbox{and}\quad \yG=y_{|\Gamma}\quad \mbox{a.e. on \,$\Sigma$},
\\
\label{sysG13}
&y(0)\,=\,0 \quad\mbox{a.e. in }\,\Omega,\quad \yG(0)=0\quad\mbox{a.e. on }\,\Gamma,
\end{align}
is obtained from \eqref{aux1}--\eqref{aux3} for $\lambda_1=\lambda_3=0, \, \lambda_2=1$.
Thanks to Lemma 3.2, it has for each ${\bf h}=(h,\hG) \in\CU$ a unique solution
$\by=(y,\yG) \in {\cal Z}$, and the associated linear mapping
\,\,${\cal G}_Q:\CU\to\CZ, \quad {\bf h}\mapsto \by$,
is continuous. The second system reads
\begin{align}
\label{sysG21}
&\dt y-\Delta y\,=\,0 \quad\mbox{a.e. in \,$Q$},\\
\label{sysG22}
&\dt\yG-\Delta_\Gamma\yG+\dn y\,=\,0\quad\mbox{a.e. on \,$\Sigma$},
\\
\label{sysG23}
&y(0)\,=\,y_0\quad\mbox{a.e. in }\,\Omega,\quad \yG(0)\,=\,y_{0_\Gamma}\quad\mbox{a.e.
on }\,\Gamma,
\end{align}
and results from \eqref{aux1}--\eqref{aux3} for $\lambda_1=\lambda_2=0,\lambda_3=1$.
For each ${\bf y}_0 \in \CN$, it has a unique solution $\by \in {\cal Z}$, and the associated mapping
\,\,${\cal G}_\Omega: {\mathcal{N}} \to {\cal Z}, \quad {\bf y}_0 \mapsto \by$,
is linear and continuous as well.
In addition, we define on the open set ${\cal A}:=({\cal U}_R\times\Phi)
\subset ({\cal U}\times\CZ)$ the nonlinear mapping
\begin{align}
\label{defG3}
&{\cal G}:{\cal A}\to \CU, \quad (\bu,\by)\mapsto {\bf h}:= (-f'(y)+u, -\fG'(\yG)+\uG)
\end{align}
as a mapping from $\CU \times \CZ$ to $\CU$.

The solution $\by=(y,\yG)$ to the nonlinear state equation \State\ is the sum
of the solution to the system \eqref{sysG11}--\eqref{sysG13}, where ${\bf h}=(h,\hG)$ is given by \eqref{defG3} 
(with $(y,\yG)$ considered as known), and of the solution to the system
\eqref{sysG21}--\eqref{sysG23}, that is, the state $\by$ associated with the control $
\bu=(u,\uG)$ is the unique solution to the nonlinear equation
\begin{equation} 
\label{nonlineq}
\by= {\cal G}_Q \big({\cal G}(\bu,\by)\bigr) +  {\cal G}_\Omega(\by_0).
\end{equation}
Let us now define  the nonlinear mapping  $\,{\cal F}:{\cal A}\to \CZ$,
\begin{align}
\label{defF}
 {\cal F}(\bu,\by)\,:=\,{\cal G}_Q\bigl({\cal G}
(\bu,\by)\bigr)+{\cal G}_\Omega(\by_0) - \by.
\end{align} 
With $ {\cal F}$, the state equation can be shortly written as
\begin{equation} \label{nonlineq2}
{\cal F}(\bu,\by)=\mathbf{0}.
\end{equation}
This equation just means that $\by=(y,\yG)$ is a solution to the state system \eqref{ss1}--\eqref{ss3} such that 
$(\bu,\by)\in{\cal A}$. From Theorem 2.2 we 
know that such a solution exists for every $\bu\in{\cal U}_R$. A fortiori, any such solution automatically enjoys the separation property 
\eqref{separation} and is uniquely determined.  

We are going to apply the implicit function theorem to the equation  \eqref{nonlineq2}. To this end,
we need the differentiability of the involved mappings.
Observe that, owing to the differentiability properties of the
involved Nemytskii operators (see, e.g., \cite[Thm.~4.22, {p.~229}]{Fredibuch}), the mapping $\,{\cal G}\,$ is twice continuously Fr\'echet differentiable in ${\cal{U}} \times \Phi$ 
as a mapping from ${\cal{U}} \times \CU$ into $\CU$,
 and for the first partial derivatives at any point 
$\,(\bus,\bys)\in {\cal A}$, and for all $\bu\in{\cal U}$
and $\by\in \CZ$, we have the identities
\begin{align}
\label{Freu}
D_{\bu}{\cal G}(\bus,\bys)[\bu]\,=\,
(u,\uG), \quad D_{\by}{\cal G}(\bus,\bys)[\by] = (-f''(\ys)y, -\fG''(\ysG)\yG). 
\end{align}  
At this point, we may  apply the chain rule, which yields  that ${\cal F}$ 
is twice continuously Fr\'echet differentiable in ${\cal U}_R\times\Phi$ as a mapping from 
$\CU\times ({\cal Y}\cap\CU)$	into $\CZ$, with the first-order partial derivatives
\begin{align}
D_{\bu}{\cal F}(\bus,\bys)\,=\,{\cal G}_Q\circ D_{\bu}{\cal 
 G}(\bus,\bys), \quad D_{\by}{\cal F}(\bus,\by)\,=\,{\cal G}_Q\circ
D_{\by}{\cal G}(\bus,\bys)-I_{\CZ},
\label{DFy}
\end{align}
where $\,I_{\CZ}\,$ denotes the identity mapping on $\,\CZ$.

We want to prove the differentiability of the control-to-state mapping ${\bf u} \mapsto  {\bf y}$  defined implicitly by the equation $\,{\cal F}({\bf u}, {\bf y})=\mathbf{0}$, using the implicit function
theorem. Now let ${\bf u^*}\in  {\cal U}_R$ be given and ${\bf y^*}={\cal S}({\bf u^*})$. We need to show that the linear
and continuous operator $\,D_{{\bf y}}{\cal F}({\bf u^*},{\bf y^*})$ is a topological isomorphism from $\CZ$ into itself.
 
To this end, let ${\bf v}\in\CZ$ be arbitrary. Then the identity $\,D_{\bf y}{\cal F}({\bf u^*},{\bf y^*})[{\bf y}]={\bf v}\,$ 
just means that $\,{\cal G}_Q\left(D_{\bf y}{\cal G}({\bf u^*},{\bf y^*})[\bf y]\right)-{\bf y}={\bf v}$, which is equivalent to saying that   
\begin{equation*}
{\bf w}\,:=\,{\bf y}+{\bf v}= {\cal G}_Q\left(D_{\bf y}{\cal G}({\bf u^*},{\bf y^*})[{\bf w}]\right)
-{\cal G}_Q\left(D_{\bf y}{\cal G}({\bf u^*},{\bf y^*})[{\bf v}]\right).
\end{equation*}
The latter identity means that ${\bf w}$ is a solution to \eqref{aux1}--\eqref{aux3} for $\lambda_1=\lambda_2=1,
\lambda_3=0$, with the specification $(h,\hG)=-D_{\bf y}{\cal G}({\bf u^*},{\bf y^*})
[{\bf v}]=(f''(\ys)v,\fG''(\ysG)\vG)\in{\cal U}$. By Lemma 3.2, such a solution ${\bf w}\in\CZ$ exists and is uniquely 
determined, which shows that  $\,D_{{\bf y}}{\cal F}({\bf u^*},{\bf y^*})$ is surjective. 
At the same time, taking ${\bf v}=\mathbf{0}$, we see that the equation $D_{{\bf y}}{\cal F}({\bf u^*},{\bf y^*})[{\bf y}]=\mathbf{0}$
means that ${\bf y}$ is the unique solution to \eqref{aux1}--\eqref{aux3} for $\lambda_1=1,\lambda_2=\lambda_3=0$. 
Obviously, ${\bf y}=\mathbf{0}$, which implies that $D_{{\bf y}}{\cal F}({\bf u^*},{\bf y^*})$ is
also injective and thus, by the open mapping principle, a topological isomorphism from $\CZ$ into itself. 

We may therefore infer from the implicit function theorem (cf., e.g., \cite[Thms. 4.7.1 and 5.4.5]{Cartan} or \cite[10.2.1]{Dieu})     
that the control-to-state mapping $\CS$ is twice continuously Fr\'echet differentiable  in ${\cal U}_R$ 
as a mapping from $\CU$ into $\CZ$. 
The explicit form of the first and second Fr\'echet derivatives is given as in Theorem 2.2:  in the case, where the directions
$(h,\hG),(k,\kG)$ belong to the space $\CU$, the corresponding solutions $(\xi,\xiG)$ and $(\eta,\etaG)$ to the
linearized system \Linear\ and to the bilinearized system \Bilinear, respectively, belong to the space $\CZ$. 
In summary, we have shown the following result.

\Bthm
\label{THM:FRECHET}
Suppose that the conditions {\bf (A1)}--{\bf (A3)} are fulfilled. Then the control-to-state operator
$\,\CS\,$ is twice continuously Fr\'echet differentiable in $\CU_R$ as a mapping from $\CU$ into $\CZ$.  
Moreover, for every $(\us,\usG)\in\CU_R$ and $(h,\hG),(k,\kG)\in\CU$ the 
functions $(\xi,\xiG)=D\CS(\us,\usG)[(h,\hG)]\in\CZ$ and $(\eta,\etaG)=D^2\CS(\us,\usG)[(h,\hG),(k,\kG)]\in\CZ$
are the unique solutions to the linearized system \Linear\ and the bilinearized system \Bilinear, respectively.
\Ethm
\Brem
It is worth noting that for the argumentation used above the actual value of the constant $R>0$ defining $\CU_R$ did not matter.
It therefore follows that $\CS$ is twice continuously Fr\'echet differentiable as a mapping from $\CU$ to $\CZ$ on the entire
space $\CU$.
\Erem

\section{The optimal control problem}
\setcounter{equation}{0}

In this section, we study the optimal control problem {\bf (CP)} with the cost functional \eqref{cost}. Besides the general
postulates {\bf (A1)}--{\bf (A3)}, we  make the following 
assumptions:
\begin{description}
\item[(A4)] \,\,The constants $\beta_i$, $i=1,2,3,4$, are nonnegative and not all zero, while $\nu,\nu_\Gamma ,\alpha,
\aG$ are positive.
\item[(A5)] \,\,The target functions satisfy $y_Q\in L^2(Q)$, $y_\Sigma\in L^2(\Sigma)$, $(y_{\Omega,T},y_{\Gamma,T})\in
\CV.$
\item[(A6)] \,\,It holds $\beta_3=\beta_4.$ 
\end{description}
\Brem
The assumptions  that $\beta_3=\beta_4$ and that $(y_{\Omega,T},y_{\Gamma,T})\in \CV$ are 
 useful in order to have regular solutions to the associated
adjoint system (see below). It is not overly restrictive in view of the continuous embedding \,$(H^1(0,T;\CH)\cap L^2(0,T;\CW\cap\CV))
\subset C^0([0,T];\CV)$\, which implies that $(y(T),\yG(T))\in\CV$. 
\Erem

The following existence result can be shown with an obvious  modification of the proof of the
corresponding theorem \cite[Thm.~3.1]{CS}. It is not restricted to functions $\,j\,$ of the special form \eqref{nondiff}.
\Bthm
Suppose that {\bf (A1)}--{\bf (A5)} are fulfilled, and suppose that $j:L^2(Q)\times L^2(\Sigma)\to\erre$ is 
convex and continuous. Then the optimal control problem {\bf (CP)} admits a solution 
$(\us,\usG)\in\Uad$.
\Ethm

In the following, we often denote by $(\us,\usG)\in\Uad$ an optimal control for {\bf (CP)} and by $(\ys,\ysG)=\CS(\us,\usG)$
the associated state. For the corresponding adjoint state system we have the following result.

\Bthm
Suppose that {\bf (A1)}--{\bf (A6)} are fulfilled and let $(u,u_\Gamma) \in \CU_R$ be a control with associated state $(y,y_\Gamma)$. 
Then the associated adjoint state system
\begin{align}
\label{adj1}
&-\dt p-\Delta p + f''(y )p =\beta_1(y -y_Q) \quad\mbox{a.e. in }\,Q,\\
\label{adj2}
&-\dt \pG -\delG\pG+\dn p + f_\Gamma''(y_\Gamma )\pG=\beta_2( y_\Gamma-y_\Sigma) \quad\mbox{and} \quad p_\Gamma=p_{|\Gamma} \quad\mbox{a.e. on }\,\Sigma,\\
\label{adj3}
&p(T)=\beta_3(y(T)-y_{\Omega,T}) \quad\mbox{a.e. in }\,\Omega, \quad \pG(T)=\beta_3(y_\Gamma(T)-y_{\Gamma,T}) \quad\mbox{a.e. on }\,\Gamma,
\end{align}
\Accorpa\Adjoint adj1  adj3
has a unique solution $(p,p_\Gamma)\in \CY$. Moreover, there is a constant $K_4>0$, which depends only on $\,R\,$ and the data, such that 
\begin{align}
\|(p,p_\Gamma)\|_{\CY} &\le K_4\left(\|y-y_Q\|_{L^2(Q)} + \|y_\Gamma-y_\Sigma\|_{L^2(\Sigma)} \nonumber \right.\\
&\left.\qquad\qquad +\|y(T)-y_{\Omega,T}\|_V +\|y_\Gamma(T)-y_{\Gamma,T}\|_{V_\Gamma} \right).\label{adj4}
\end{align}
\Ethm
\Bdim With the exception of \eqref{adj4}, the assertion follows from \cite[Thm.~3.4]{CS}. To show \eqref{adj4}, we argue as follows:
we  put 
\begin{align}
&q(x,t):=p(x,T-t),\quad q_\Gamma:= p_{\Gamma}(x,T-t), \quad \tilde y(x,t):=y(x,T-t),\nonumber \\
&\tilde y_{\Gamma}(x,t):=y_\Gamma(x,T-t), \quad h(x,t):=\beta_1(y(x,T-t)-y_Q(x,T-t)),\nonumber\\
&\hG(x,t):=\beta_2(y_\Gamma(x,T-t)-y_\Sigma(x,T-t))\,.
\label{shifted}
\end{align}
In terms of these quantities, the adjoint system \Adjoint\ takes the form
\begin{align}
\label{anna1}
&\dt q-\Delta q= -f''(\tilde y)q+h \quad\mbox{a.e. in }\,Q,\\
\label{anna2}
&\dt q_\Gamma-\Delta_\Gamma q_\Gamma+\dn q=-\fG''(\tilde y_\Gamma)+\hG \quad\mbox{and}\quad q_\Gamma=q_{|\Gamma} 
\quad\mbox{a.e. on }\,\Sigma,\\
\label{anna3}
&q(0)=\beta_3(y(T)-y_{\Omega,T}) \quad\mbox{a.e. in }\,\Omega,\quad q_\Gamma(0)=\beta_3(y_\Gamma(T)-y_{\Gamma,T}) \quad\mbox{a.e. on }\,\Gamma,
\end{align}
which is a special case of the auxiliary system \eqref{aux1}--\eqref{aux3} with $\lambda_1=\lambda_2=\lambda_3=1$, up to the
minor difference that the arguments of the functions $f''$ and $\fG''$ differ from those in \eqref{aux1} and \eqref{aux2} 
by a time shift. Obviously, $(h,\hG)\in L^2(0,T;\CH)$, and it follows from the continuity of the embedding
${\cal Y}\subset C^0([0,T];\CV)$ and {\bf (A5)} that the initial data of $(q,q_\Gamma)$ belong to ${\cal V}$. Therefore, if we put  
\[
M:= \lambda_2\,\|(h,\hG)\|_{L^2(0,T;\CH)}\,+\,\lambda_3\,\|(y_0,y_{0_\Gamma}\|_{\CV},
\]
then the estimates \eqref{uwe1}--\eqref{uwe4} performed in the proof of Lemma 3.2 can be repeated (notice that the 
condition $(q(0),q_\Gamma(0))\in\CU$ is not needed for these estimates), and \eqref{adj4} follows from \eqref{uwe4}. 
\Edim
\noindent Note that -- at this point -- we cannot expect to have the stronger regularity 
$(p,p_\Gamma) \in\CZ$, since, in view of assumption {\bf (A5)},  the \rhs s of \eqref{adj1} and 
\eqref{adj2} only belong to $L^2$--spaces.

\subsection{First-order necessary optimality conditions}

In this section, we aim at deriving  associated first-order necessary optimality conditions for local minima of the
optimal control problem {\bf (CP)}. We assume that {\bf (A1)}--{\bf (A6)} are fulfilled and that $j:L^2(0,T;\CH)\to\erre$ is a general convex and continuous functional. We use the abbreviations
$${\bf u}:=(u,\uG),\quad {\bf u^*}:=(\us,\usG),\quad {\bf y}:=(y,\yG),\quad {\bf y^*}=(\ys,\ysG).$$ 
Next, we define  the reduced cost functionals associated with the functionals  $J$ and ${\cal J}$  introduced in \eqref{diffpart} and 
\eqref{cost} by
\begin{equation}\label{reduced}
\widehat J(\bu) =  J(\CS(\bu),{\bf u}), \quad \wJ(\bu)={\cal J}(\CS(\bu),\bu)\,.
\end{equation}
Since $\CS$ is twice continuously Fr\'echet differentiable  from $\CU$ into the space $C^0([0,T];\CH)\subset\CZ$, 
it follows from the chain rule that $\,\widehat J\,$ is  a twice continuously Fr\'echet 
differentiable mapping from $\CU$ into $\erre$, where, for every $\bus=(\us,\usG)\in\CU$ and every ${\bf h}=(h,\hG)\in\CU$, it holds 
with $(\ys,\ysG)=\CS(u^*,\usG)$ that
\begin{align}
\label{DJ}
D\widehat J(\bus)[{\bf h}]\,&=\,\beta_1\iint_Q \xi(y^*-y_Q)\,+\,\beta_2\iint_\Sigma \xi_{\Gamma}(\ysG-y_\Sigma)\,+\,\beta_3\iO
\xi(T)(\ys(T)-y_{\Omega,T})  \nonumber\\
&\quad + \beta_3\int_\Gamma \xi_{\Gamma}(T)(\ysG(T)-y_{\Gamma,T})\,+\,\nu \iint_Q \us h \,+\,\nu_\Gamma \iint_\Sigma \usG \hG,
\end{align}
where $\,(\xi,\xi_{\Gamma})=D\CS(\bus)[\bf h]\in\CY\,$ is the unique solution to the linearized system \Linear\
associated with ${\bf h}$.
\Brem
Observe that the \rhs\ of \eqref{DJ} is meaningful also for arguments ${\bf h}=(h,\hG)\in L^2(0,T;\CH)$, where in this case $(\xi,\xi_\Gamma)
=D\CS(\bus)[{\bf h}]$ with the extension of the operator $D\CS(\bus)$ to $L^2(0,T;\CH)$ introduced in Remark 3.1. Hence, by means of the
identity \eqref{DJ} we can extend the operator $D\widehat J(\bus)\in\CU^*$ to $L^2(0,T;\CH)$. The extended operator, which we again denote by
$D\widehat J(\bus)$, then becomes an element of $(L^2(0,T;\CH))^*$. In this way, expressions of the form $D\widehat J(\bus)[{\bf h}]$  have a proper meaning also for 
${\bf h}\in L^2(0,T;\CH)$.  
\Erem
In the following, we assume that $\bus = (\us,\usG)$ is a given locally optimal control for ${\bf (CP)}$ in the sense
of \,$\CU$, that is, there is some $\varepsilon>0$ such that
\begin{equation}
\label{lomin}
\wJ(\bu)\,\ge\,\wJ(\bus)\quad\mbox{for all $\bu\in\Uad$ satisfying }\,\|\bu-\bus\|_{\CU}\,\le\,\varepsilon.
\end{equation}

Notice that any locally optimal control in the sense of $L^p(Q)\times L^p(\Sigma)$ with $1 \le p < \infty$
is also locally optimal in the sense of $\CU$, since the topology of $\CU$ is the finest among these spaces. 
Therefore, a result proved for locally optimal controls in the sense of  $\CU$ is also valid for locally optimal controls 
in the sense of  $L^p(Q)\times L^p(\Sigma)$. It is also true for  (globally) optimal controls.

We claim that then  the variational inequality
\begin{equation}
\label{var1}
{D{\widehat J}(\bus)[\bu -  \bus] + j(\bu) - j(\bus) \ge 0 \quad \forall\, \bu \in \Uad}
\end{equation}
is satisfied. Although this result seems to be pretty standard by now, we nevertheless provide its proof for the reader's convenience.
To this end, note that by local optimality there is some $\varepsilon>0$ such that \eqref{lomin} is satisfied.

Now let $\bu\in\Uad$ be arbitrary. Then, for any $\tau\in (0,1]$,  we have $\bu_\tau:=\bus+\tau(\bu-\bus)\in\Uad$. For all sufficiently small $\tau > 0$, it holds in addition that 
$\|\bu_\tau-\bus\|_{\CU}\,\le \varepsilon.$
Hence, for all such sufficiently small $\tau>0$, we obtain from \eqref{lomin} and from the convexity of $j$ the following chain of inequalities:
\begin{align*}
0\,&\le \,\wJ(\bu_\tau)-\wJ(\bus)\,\le\,\widehat J(\bu_\tau)-\widehat J(\bus)+j(\bu_\tau)-j(\bus)\\
&\le J(\CS(\bu_\tau),\bu_\tau)-J(\CS(\bus),\bus) + \tau(j(\bu)-j(\bus))\,.
\end{align*}
Now, we divide by $\tau>0$ and take the limit as $\tau\searrow0$. Since $\,\|\bu_\tau-\bus\|_\CU\to 0\,$ 
as $\,\tau\searrow0$, we then 
can infer that the claim \eqref{var1} is actually valid.

The variational inequality \eqref{var1}, in turn, implies that \,$\bus$\, solves the convex minimization problem
\[
\min_{\bu \in \CU} \,\bigl(\Psi(\bu)+ j(\bu) + I_{\Uad}(\bu)\bigr),
\]
with $\,\Psi(\bu)=D\widehat J(\bus)[\bu]$, and where $I_{\Uad}$ denotes the indicator function of $\Uad$. Hence,
denoting by the symbol $\,\partial\,$ the subdifferential mapping in $L^2(0,T;\CH)$ (recall that $j$ is a convex 
continuous functional on  $L^2(0,T;\CH)$), we have the inclusion 
$\mathbf{0}\in\partial \bigl(\Psi+j+I_{\Uad}\bigr)(\bus)$  or, by the well-known rules for 
subdifferentials of convex functionals, 
$$ \mathbf{0}\in \{D{\widehat J}(\bus)\}+\partial j(\bus)+\partial I_{\Uad}(\bus).$$
In other words, there
are ${\boldsymbol \lambda}^*\in\partial j(\bus)$ and $\widehat{\boldsymbol \lambda}\in\partial I_{\Uad}(\bus)$ such that
$\,\mathbf{0}=D{\widehat J}(\bus)+{\boldsymbol \lambda}^*+\widehat{\boldsymbol \lambda}$, which by the definition
of $\partial I_{\Uad}(\bus)$ means that
\[
0\le D{\widehat J}(\bus)[\bu-\bus]+{\boldsymbol \lambda}^*[\bu-\bus]\quad\forall\,\bu\in\Uad.
\] 
We have thus shown the following result (where we identify $\bl^*$  
with the corresponding element of $\,L^2(0,T;\CH)\,$ 
according to the Riesz isomorphism): if $\bus\in\Uad$ is a locally optimal control for {\bf (CP)} in the sense of 
$\CU$, then there is 
some $\bl^*=(\lambda^*,\lambda_\Gamma^*) \in \partial j(\bus)\subset L^2(0,T;\CH)$ such that 
\begin{equation} \label{varineq1}
D{\widehat J}(\bus)[\bu - \bus] +  \iint_Q 
\lambda^* (u- \us)\,+\iint_\Sigma \lambda_\Gamma^*(\uG-\usG) \,\ge 0 \quad \forall \,
\bu=(u,\uG) \in \Uad.
\end{equation}
As usual, we simplify the expression $D{\widehat J}(\bus)[\bu-\bus]$ in \eqref{varineq1}
by means of the adjoint state variables defined in \Adjoint. A standard calculation (see \cite[Thm.~3.4]{CS}) then leads 
to the following result.
\Bthm
\,\,{\rm (Necessary optimality condition)}  \label{Thm4.5}Suppose that {\bf (A1)}--{\bf (A6)}
are fulfilled and that $j:L^2(0,T;\CH)\to\erre$ is convex and continuous. Moreover, 
 let $\bus=(\us,\usG) \in \Uad$ be a locally optimal control of  {\bf (CP)} in the sense of $\,\CU\,$ 
with associated state $\,(\ys,\ysG)={\cal S}(\bus)$
and adjoint state $\,{\bf p^*}=(\ps,\psG)$.
Then there exists some ${\boldsymbol \lambda}^*=(\lambda^*,\lambda_\Gamma^*)   \in \partial j(\bus)$ such that,
for all $\, \bu =(u,\uG)\in \Uad$, 
\begin{align}
\label{varineq2}
&\iint_Q  \left(p^*+\lambda^*+\nu u^*\right)\left(u-\us\right) \,+\,\iint_\Sigma\left(p_\Gamma^*+\lambda_\Gamma^*+\nu_\Gamma  \usG\right)
\left(\uG-\usG\right) \, \ge \,0 \,.
\end{align}
\Ethm
We underline again that \eqref{varineq2} is also necessary for all globally optimal controls and all controls 
which are locally optimal in the sense of $L^p(Q)\times L^p(\Sigma)$ with $p \ge 1$.

\subsection{Sparsity of controls}

The convex function $j$ in the objective functional accounts for the sparsity of optimal controls, i.e., 
any locally optimal control can vanish in some region of the space-time cylinder $Q$. The form of this region depends on the particular choice of the functional $j$ which can differ in different situations.
The sparsity properties can be deduced from the variational inequality \eqref{varineq2} and the particular
form of the subdifferential  $\partial j$. In this paper, we are mostly interested in the so-called 
{\em directional sparsity with respect to time} and {\em full sparsity}, and we omit the case of {\em directional
sparsity with respect to space} which can be handled analogously. Our analysis closely follows the lines of 
\cite[Sect.~4]{SpTr1}. We therefore can be brief. 

In order to have directional sparsity with respect to time, we use the functionals
\begin{align} 
&j_T^Q: L^1(0,T;L^2(\Omega)) \to \mathbb{R},\quad j^\Sigma_T:L^1(0,T;L^2(\Gamma)) \to \erre,\nonumber\\[1mm]
\label{jT}
&j_T^Q(u) = \int_0^T\,\|u(\cdot,t)\|_{\Ldue}\,dt, \quad j_T^\Sigma = \int_0^T \|\uG(\cdot,t)\|_{L^2(\Gamma)}\,dt.
\end{align}
The associated subdifferential of $\,j_T^Q\,$ is given by (cf., \cite{herzog_stadler_wachsmuth2012})
\begin{equation}
\label{djTQ}
\partial j_T^Q(u) = \left\{\lambda \in L^2(Q): 
\left\{
\begin{array}{ll}
\|\lambda(\cdot,t)\|_{\Ldue} \le 1 &\,\mbox{ if } u(\cdot,t) = 0\\[1ex]
\lambda(\cdot,t)\,=u(\cdot,t)/ \|u(\cdot,t)\|_{\Ldue}&\,\mbox{ if } u(\cdot,t) \not= 0
\end{array}
\right\} 
\right\}\,,
\end{equation}
where the properties above are satisfied for a.e. $t \in (0,T)$. The subdifferential of $\,j_T^\Sigma\,$ is obtained analogously.

The case of full sparsity is obtained for the functionals 
\begin{align} \label{g3}
&j_F^Q: L^1(Q) \to \mathbb{R}, \quad j_F^\Sigma:L^1(\Sigma)\to\erre, \nonumber\\[1mm] 
&j_F^Q(u) = \|u\|_{L^1(Q)}, \quad j_F^\Sigma(u)=\|u\|_{L^1(\Sigma)}.
\end{align}
In this case we have (see \cite{ioffe_tikhomirov1979})
\begin{equation}
\label{djF}
\partial j_F^Q(u) = \left\{\lambda \in L^2(Q):\,
\lambda(x,t) \in \left\{
\begin{array}{ll}
\{1\} & \mbox{ if } u(x,t) > 0\\
{[-1,1]}& \mbox{ if } u(x,t) = 0\\
\{-1\} & \mbox{ if } u(x,t) < 0\\
\end{array}
\right\}
\mbox{ for a.e. } {(x,t) \in Q}
\right\}\,,
\end{equation}
and the subdifferential of $\,j_F^\Sigma\,$ is obtained analogously.

\subsubsection{Directional sparsity in time} 

In this section, we will focus on directional sparsity in time. 
To this end, we discuss the following auxiliary variational inequality:
\begin{equation} \label{varineqaux}
\iint_Q (d(x,t) + \kappa \lambda(x,t) + \nu u(x,t))(v(x,t) - u(x,t))\,dx\,dt \ge 0 \quad \forall\, v \in C,
\end{equation}
where $\lambda \in \partial j_T^Q(u)$ and 
\begin{equation} \label{C}
C = \{v \in L^\infty(Q): \underline u \le v(x,t) \le {\overline u} \mbox{ a.e. in } Q\}
\end{equation}
with given real numbers $\underline u < 0 < \overline u$, $\kappa> 0$, $\nu > 0$, and a given function
$d \in L^2(Q)$. 

The following result is known from \cite{herzog_stadler_wachsmuth2012,casas_herzog_wachsmuth2017,SpTr1}.
\Blem
\label{L4.1} \,\,{\rm (Sparsity)} \quad Let $\underline u < 0 < \overline u$, $\kappa> 0$, $\nu > 0$, and let $\,u \in C\,$ be a solution {to} the variational 
inequality \eqref{varineqaux}. Then, for a.e. $t\in (0,T)$,
\begin{equation}\label{equivalence}
u(\cdot,t) = 0 \quad \Longleftrightarrow \quad  \|d(\cdot,t)\|_{L^2(\Omega)} \le \kappa,
\end{equation}
as well as 
\begin{equation} \label{subdiff}
\lambda(\cdot,t) \left\{
\begin{array}{lcl}
\in \{v\in\Ldue:\|v\|_{\Ldue}\le 1\}&\,\mbox{ if } \, \|u(\cdot,t) \|_{L^2(\Omega)} = 0\\[1ex]
= \displaystyle \frac{u(\cdot,t)}{\|u(\cdot,t) \|_{L^2(\Omega)}}&\,\mbox{ if } \, \|u(\cdot,t) \|_{L^2(\Omega)} \not= 0
\end{array}
\right. \,.
\end{equation}
\Elem
\Brem
A corresponding sparsity result can be obtained for functions defined on the lateral boundary $\Sigma$ if the variational inequality 
\eqref{varineqaux} and the set $\,C\,$ are adapted accordingly. For the sake of a shorter exposition, and since the necessary changes are obvious, 
we may leave it to the reader to formulate the details.   
\Erem
We apply the lemma, and its counterpart for functions defined on $\Sigma$, to derive sparsity properties of locally  optimal 
controls from the variational inequality \eqref{varineq2}. For directional sparsity in time, we use the convex and continuous functional
\begin{equation} \label{gT}
j(\bu) =  j((u,\uG)) : = \alpha\,j_T^Q(u) + \aG\,j_T^\Sigma(\uG) = \alpha\,j_T^Q(I(\bu)) + \aG\,j_T^\Sigma(I_\Gamma(\bu)),
\end{equation}
where $I,I_\Gamma$ denote the linear and continuous projection mappings $I: \bu=(u,\uG)  \mapsto u$ and $I_\Gamma:\bu=(u,\uG)\mapsto \uG$ from
$L^2(Q)\times L^2(\Sigma)$ to $L^2(Q)$ and $L^2(\Sigma)$, respectively.
Since the convex functionals $j_T^Q$ and $j_T^\Sigma$ are continuous on the whole spaces $L^2(Q)$ and $L^2(\Sigma)$, respectively,
we obtain from the rules for subdifferentials (cf., 
\cite[Sect.~4.2.2,~Thms.~1~and~2]{ioffe_tikhomirov1979}) that 
\begin{align*}
\partial j(\bu)\, &=\,\alpha\, I^*\, \partial j_T^Q(I(\bu))\, +\,\aG\, I_\Gamma^*\, \partial j_T^\Sigma(I_\Gamma(\bu)) \\
&= \left\{ (\alpha\lambda,\aG\lambda_\Gamma) \in L^2(Q)\times L^2(\Sigma): \lambda \in \partial j_T^Q(u), 
\quad \lambda_\Gamma\in\partial j_T^\Sigma(\uG)\right\}.
\end{align*}

The variational inequality \eqref{varineq2} is equivalent to two independent variational inequalities for $\us$ and $\usG$ 
that have to hold simultaneously, namely,
\begin{eqnarray}
\iint_Q  \left( p^*+\alpha\lambda^*+\nu \us\right)\left(u - \us\right) &\!\!\ge\!\!& 0 \quad \forall\, u \in U^{\rm ad}, \label{varin1}\\
\iint_\Sigma  \left(p^*_\Gamma+\aG\lambda^*_\Gamma +\nu_\Gamma \usG\right)\left(\uG - \usG\right) 
&\!\!\ge\!\!& 0 \quad \forall\, \uG \in U^{\rm ad}_\Gamma, \label{varin2}
\end{eqnarray}
where 
\begin{align}
\label{uaduno}
U^{\rm ad}\,& = \,\{u \in L^\infty(Q): \rhomin  \le u(x,t) \le \,\rhomax \mbox{ for a.e. } (x,t) \in Q\},\\
\label{uaddue}
U^{\rm ad}_\Gamma\,&=\,\{\uG\in \LiS: \Rmin \le \uG(x,t) \le\Rmax \,\mbox{ for a.e. }\,(x,t)\in \Sigma\},
\end{align} 
and where, for a.e. $t\in\oT$, 
\begin{equation} 
\label{subdi1}
\lambda^*(\cdot,t) \left\{
\begin{array}{lcl}
\in \{v\in \Ldue: \,\|v\|_{\Ldue}\,\le\,1\}&\,\mbox{ if } \, \|\us(\cdot,t) \|_{L^2(\Omega)} = 0\\[1ex]
= \,\displaystyle \frac{\us(\cdot,t)}{\|\us(\cdot,t) \|_{L^2(\Omega)}}&\,\mbox{ if } \, \|\us(\cdot,t) \|_{L^2(\Omega)} \not= 0
\end{array}
\right. \,.
\end{equation}
as well as
\begin{equation} 
\label{subdi2}
\lambda^*_\Gamma(\cdot,t) \left\{
\begin{array}{lcl}
\in \{\vG\in L^2(\Gamma): \|\vG\|_{L^2(\Gamma)}\,\le\,1\}&\,\mbox{ if } \, \|\usG(\cdot,t) \|_{L^2(\Gamma)} = 0\\[1ex]
= \,\displaystyle \frac{\usG(\cdot,t)}{\|\usG(\cdot,t) \|_{L^2(\Gamma)}}&\,\mbox{ if } \, \|\usG(\cdot,t) \|_{L^2(\Gamma)} \not= 0
\end{array}
\right. .
\end{equation}
Applying Lemma \ref{L4.1} to \eqref{varin1}, and its analogue on $\Sigma$ to \eqref{varin2}, we arrive at the following result:
\Bthm  {\rm (Directional sparsity in time)} \,\,\, Suppose that the general assumptions {\bf (A1)}--{\bf (A6)} are fulfilled, and assume 
that \,$\rhomin<0<\rhomax\,$ and $\,\Rmin<0<\Rmax$. Let $\bus= (\us,\us_\Gamma)\in\Uad$ be 
a locally optimal control in the sense of \,${\cal U}$\, of the problem {\bf (CP)}  
with the sparsity functional $\,j\,$ defined in \eqref{jT}, and with associated state $(\ys,\ysG)=\CS(\bus)$
solving \State\ and adjoint state ${\bf p}^*= (p^*,\psG)$ solving  \eqref{adj1}--\eqref{adj3}. 
Then there are functions $\lambda^*,\lambda^*_\Gamma$ that satisfy \eqref{subdi1}, \eqref{subdi2}, and \eqref{varin1}--\eqref{varin2}.
In addition, for
almost every $t\in (0,T)$, we have that
\begin{eqnarray}
\|\us(\cdot,t)\|_{L^2(\Omega)} = 0 \quad &\Longleftrightarrow& \quad \|p^*(\cdot,t)\|_{L^2(\Omega)} \le \alpha, 
\label{usparsity}
\\
\label{uGsparsity}
\|\usG(\cdot,t)\|_{L^2(\Gamma)} = 0 \quad &\Longleftrightarrow& \quad \|\psG(\cdot,t)\|_{L^2(\Gamma)} \le \aG. 
\end{eqnarray}
Moreover, if\, ${\bf p}^*$ and $\lambda^*, \lambda^*_\Gamma$ are given, then
the optimal controls
$\us$,  $\usG$ are obtained from the projection formulas
\begin{eqnarray*}
\us(x,t)& =& \max\left\{\rhomin, \min\left\{ \rhomax, -{\nu}^{-1} \left(p^*+ \alpha \,\lambda^*\right)(x,t)\right\}\right\}
\,\mbox{ for a.e. $(x,t)\in Q$} ,\\
\usG(x,t) &=& \max\left\{\Rmin,\min\left\{ \Rmax,-{\nu_\Gamma }^{-1} \left(\psG+\aG\,\lambda_\Gamma^*\right)(x,t)\right\}\right\}
\,\mbox{ for a.e. $(x,t) \in \Sigma$.}
\end{eqnarray*}
\Ethm

The projection formulas above are standard conclusions from the variational inequalities \eqref{varin1}--\eqref{varin2}.
It is to be expected that the support of locally optimal controls $(\us,\usG)$ will shrink with increasing sparsity   parameters $\alpha,\aG$.
Although this can hardly be quantified or proved, it is useful to confirm that optimal controls vanish for all sufficiently large values of $\alpha$
and $\aG$. We are going to derive a corresponding result now.
 
For this purpose, let us indicate for a while the dependence of optimal controls, optimal states, and the associated
adjoint states, on the pair of weights ${\boldsymbol \alpha}:=(\alpha,\aG)$ by an index $\bal$, i.e., we write
$$\bu^*_\bal=(\us_\bal, \us_{\Gamma_\bal}), \quad {\bf y}^*_\bal=(\ys_\bal,\ys_{\Gamma_\bal})\quad\mbox{and}\quad 
{\bf p}^*_\bal=(p^*_\bal,\ps_{\Gamma_\bal}).$$
From \eqref{usparsity} 
and \eqref{uGsparsity} we infer that 
$\us_\bal = 0 $ holds for all
$\alpha > \alpha^*$ if 
\begin{equation} 
\label{globalbound1}
\alpha^* := \sup_{\alpha> 0} \| p^*_\bal\|_{L^\infty(0,T;H)} < +\infty,
\end{equation}
and $\us_{\Gamma_\bal} = 0 $ holds for all
$\aG > \alpha^*_\Gamma$ if 
\begin{equation} 
\label{globalbound2}
\alpha_\Gamma^* := \sup_{\aG > 0}  \|\ps_{\Gamma_\bal}\|_{L^\infty(0,T;\HG)} < +\infty.
\end{equation}
Hence the controls vanish simultaneously if the components of $\bal=(\alpha,\aG)$ are both larger than the corresponding components 
of $\bal^*:=(\alpha^*,\alpha^*_\Gamma)$ provided that both \eqref{globalbound1} and \eqref{globalbound2} are valid. 
We now show that this is actually the case. 
To this end, we recall the global
estimates \eqref{ssbound1}--\eqref{ssbound2}, which have to be satisfied by 
all possible  states corresponding to controls $\bu\in\Uad$. Consequently, the \rhs s 
$h:=\beta_1 (\ys - y_Q)$ and $\hG:=\beta_2 (\ysG - y_\Sigma)$ in \eqref{adj1} and \eqref{adj2} 
are uniformly bounded in $L^2(Q)$ and $L^2(\Sigma)$, respectively,
independently of $\bal$. Moreover, owing to the continuity of the embedding $\,{\cal Y}\subset C^0([0,T];\CV)$, the
terminal data $(p^*(T),\psG(T))=(\beta_3(\ys(T)-y_{\Omega,T}),\beta_3(\ysG(T)-y_{\Gamma,T}))$ are uniformly bounded in $\CV$.
Therefore, it follows from \eqref{adj4} that $(p^*,\psG)$ is bounded in ${\cal Y}$, independently of $\bal$. In particular,
$$
\|(p^*,\psG)\|_{C^0([0,T];\CH)}\,\le\,C\,,
$$
where $C>0$ is independent of $\bal$. 
Thus, in the case of directional sparsity in time, locally optimal controls in the sense of 
\,$\CU$\, vanish for
 sufficiently large sparsity parameters.

\subsubsection{Full sparsity}

In this section, we consider the case when the sparsity functional is given by \eqref{nondiff}, i.e.,
\Beq
\label{elvis}
j(\bu)=j((u,\uG)):=\alpha \,j_F^Q(u)+\aG\,j_F^\Sigma(\uG)=\alpha\,j_F^Q(I(\bu))+\aG\,j_F^\Sigma(I_\Gamma(\bu)),
\Eeq
where $I$ and $I_\Gamma$ have the same meaning as in \eqref{gT}. Similarly as there, we obtain
that  
$$
\partial j(\bu)=\{(\alpha\lambda,\aG\lambda_\Gamma):\lambda\in\partial j_F^Q(u),\quad \lambda_\Gamma\in\partial j_F^\Sigma(\uG)\}.
$$
In the problem of full sparsity, the variational inequality \eqref{varineq2} becomes
\begin{align} \label {varineq_full}
&\iint_Q  \left(p^*+\alpha \lambda^*+\nu u^*\right)\left(u-\us\right) \,+\,\iint_\Sigma\left(p_\Gamma^*+\alpha_\Gamma \lambda_\Gamma^*+\nu_\Gamma  \usG\right)
\left(\uG-\usG\right) \, \ge \,0 \nonumber\\
&\qquad\mbox{for all }\,(u,u_\Gamma) \in \Uad.
\end{align}
We now show the following result.
\Bthm {\rm (Full sparsity)} \,\,\,Suppose that the assumptions {\bf (A1)}--{\bf (A6)} are fulfilled, and assume that 
$\rhomax <0 <\rhomin$ and $\Rmin <0 <\Rmax$. Let $\bus=(\us,\usG)\in \Uad$ be a locally optimal control in the sense 
of \,$\CU$\, for the problem {\bf (CP)}
with the cost functional $\,j\,$ defined in \eqref{elvis}, and with associated state $(\ys,\ysG)=\CS(\bus)$ solving \State\ and 
adjoint state ${\bf p}^*=(p^*,\psG)$ solving \Adjoint. Then there exist
functions $\lambda^*\in\partial j_F^Q(\us)$ and $\lambda^*_\Gamma\in \partial j_F^\Sigma(\usG)$ that satisfy \eqref{varin1}--\eqref{varin2}.
In addition, we have that
\begin{eqnarray}
\us(x,t) = 0 \quad &\Longleftrightarrow& \quad |p^*(x,t)| \le \alpha, \quad\mbox{for a.e. }\,(x,t)\in Q, 
\label{ufullsparsity}
\\
\label{uGfullsparsity}
\usG(x,t) = 0 \quad &\Longleftrightarrow& \quad |\psG(x,t)| \le \aG, \quad\mbox{for a.e. }\,(x,t)\in \Sigma. 
\end{eqnarray}
Moreover, if\, ${\bf p}^*$ and $\lambda^*, \lambda^*_\Gamma$ are given, then
the optimal controls
$\us$,  $\usG$ are obtained from the projection formulas
\begin{eqnarray*}
\us(x,t)& =& \max\left\{\rhomin, \min\left\{ \rhomax, -{\nu}^{-1} \left(p^*+ \alpha \,\lambda^*\right)(x,t)\right\}\right\}
\,\mbox{ for a.e. $(x,t)\in Q$} ,\\
\usG(x,t) &=& \max\left\{\Rmin,\min\left\{ \Rmax,-{\nu_\Gamma }^{-1} \left(\psG+\aG\,\lambda_\Gamma^*\right)(x,t)\right\}\right\}
\,\mbox{ for a.e. $(x,t) \in \Sigma$.}
\end{eqnarray*}
\Ethm
\Bdim
First, we observe that the projection formulas are a direct consequence of the variational inequalities \eqref{varin1} and \eqref{varin2}.
It thus only remains to show the validity of \eqref{ufullsparsity} and \eqref{uGfullsparsity}. We only prove the former equivalence, the
proof of the latter is analogous.  

We use the first projection formula and the fact that $\rhomin<0<\rhomax$. For a.e. $(x,t)\in Q$, we have: if $\,\us(x,t)=0$, then
$\,\,-\nu^{-1}(p^*(x,t)+\alpha\lambda^*(x,t))=0$, where $\lambda^*(x,t)\in [-1,1]$. Consequently,
$\,|p^*(x,t)|=\alpha|\lambda^*(x,t)|\le\alpha$.

Now let us assume that $|p^*(x,t)|\le \alpha$. If $\us(x,t)>0$, then $\lambda^*(x,t)=1$ and
$\,\,-\nu^{-1}(p^*(x,t)+\alpha)\ge u^*(x,t)>0,$ which implies that $\,p^*(x,t)+\alpha<0\,$ and thus $\,|p^*(x,t)|=-p^*(x,t)>\alpha$,
a contradiction. By analogous reasoning, we can show that also the assumption $\,u^*(x,t)<0$ \,leads to a contradiction. 
We thus must have $u^*(x,t)=0$. 
This ends the proof.
\Edim

We conclude this section by investigating whether optimal controls have to vanish for sufficiently large sparsity parameters. With the denotation
introduced in the previous section, we thus have to check whether 
\begin{eqnarray}
\label{presley}
\alpha^*:=\sup_{\alpha>0}\,\|p^*_\bal\|_{\LiQ}\,<+\infty,\quad \alpha^*_\Gamma:=\sup_{\alpha_\Gamma>0}\,\|p^*_{\Gamma_\bal}\|_{\LiS}
\,<+\infty\,.
\end{eqnarray}
Such  bounds cannot be expected to hold, in general. But they are actually valid under the following additional assumption:

\begin{description}
\item[(A7)] \,It holds $\beta_3=\beta_4=0$, as well as $y_Q\in\LiQ$ and $y_\Sigma\in \LiS$.
\end{description}

\noindent Indeed, if {\bf (A7)} is fulfilled, then the quantities introduced in \eqref{shifted} satisfy \eqref{anna1},
\eqref{anna2}, as well as $q(0)=0$ and $q_\Gamma(0)=0$, where
the functions $\,h\,$ and $\,\hG\,$ are bounded in $\LiQ$ and $\LiS$, respectively, independently of $\bal$.
Now observe that in terms of these quantities the adjoint system \Adjoint\ becomes
a special case of the auxiliary system \eqref{aux1}--\eqref{aux3} with $\lambda_1=\lambda_2=1$ and $\lambda_3=0$, up to the
minor difference that the arguments of the functions $f''$ and $\fG''$ differ from those in \eqref{aux1} and \eqref{aux2} 
by a time shift. Since this difference 
does not matter in the estimates performed 
in the proof of Lemma 3.2, we may argue as there to conclude that
$$\|(p^*_\bal,p^*_{\Gamma_\bal})\|_\CU \,=\,\|(q,q_\Gamma)\|_\CU \,\le\,C_1\,\|(h,\hG)\|_\CU\,\le\,C_2,$$
where $C_1$ and $C_2$ do not depend on $\bal$. The condition \eqref{presley} is therefore fulfilled.
In conclusion, also in this case all locally optimal controls in the sense of \,$\CU$\, vanish for sufficiently large sparsity parameters.

\subsection{Second-order sufficient optimality conditions} 
We conclude this paper with the derivation of second-order sufficient optimality conditions. 
We provide conditions that ensure local optimality of pairs \,$\bus=(u^*,u^*_\Gamma)$ obeying the first-order necessary optimality 
conditions of Theorem \ref{Thm4.5}. Second-order sufficient optimality conditions are based on a condition of coercivity that 
is required to hold for the smooth part  $\,J\,$ 
of $\,{\cal J}\,$ in a certain critical cone. The nonsmooth part $\,j\,$ contributes to sufficiency by its convexity. In the following,
we restrict ourselves to the case of full sparsity, where we generally assume that {\bf (A1)}--{\bf (A6)} and the conditions
$\,\rhomin<0<\rhomax\,$ and $\,\Rmin<0<\Rmax\,$ are fulfilled.
Our analysis will follow closely the lines of \cite{casas_ryll_troeltzsch2015}, where a second-order analysis was 
performed for sparse control of the FitzHugh--Nagumo system. In particular, we adapt the proof of 
\cite[Thm.~3.4]{casas_ryll_troeltzsch2015} to our setting of less regularity.

To this end, we fix  a pair of controls $\,\bus = (u^*,u^*_\Gamma)$ that satisfies the first-order necessary optimality conditions,
and we set $\,{\bf y}^*=(\ys,\ysG)=\CS(\bus)$. Then the cone
\[
C(\bus) = \{ (v,v_\Gamma) \in L^2(0,T;\CH) \,\text{ satisfying the sign conditions \eqref{sign} a.e. in $Q$ and $\Sigma$}\},
\]
where
\begin{equation} \label{sign}
v(x,t) \left\{
\begin{array}{l}
\ge 0 \,\, \text{ if }\,\, u^*(x,t) = \rhomin\\ 
\le 0 \,\, \text{ if }\,\, u^*(x,t) = \rhomax 
\end{array}
\right. \,,\quad 
v_\Gamma(x,t) \left\{
\begin{array}{l}
\ge 0 \,\, \text{ if }\,\, u_\Gamma^*(x,t) = \Rmin\\ 
\le 0 \,\, \text{ if } \,\,u^*(x,t) = \Rmax 
\end{array}
 \right. \,, 
\end{equation}
is called the {\em cone of feasible directions}, which is a convex and closed subset of $L^2(0,T;\CH)$.
 We also need the directional derivative of $j$ at $\bu\in L^2(0,T;\CH)$ in the direction 
$\bv\in L^2(0,T;\CH)$, which is given by
\begin{equation}
\label{j'}
j'(\bu,\bv) = \lim_{\tau \searrow 0} \frac{1}{ \tau}(j(\bu+\tau \bv)-j(\bu))\,.
\end{equation}
Following the definition of the critical cone in \cite[Sect. 3.1]{casas_ryll_troeltzsch2015}, we define
\begin{equation}
\label{critcone}
C_{\bus} = \{\bv \in C(\bus): D\widehat{J}(\bus)[\bv] + j'(\bus,\bv) = 0\}\,,
\end{equation}
which is also a closed and convex subset of $L^2(0,T;\CH)$. According to \cite[Sect. 3.1]{casas_ryll_troeltzsch2015}, 
it consists of all $\bv=(v,\vG)\in C(\bus)$ satisfying  
\begin{equation} \label{pointwise}
v(x,t) \left\{
\begin{array}{l}
= 0 \,\,\mbox{ if } \,\,|p^*(x,t) + \nu u^*(x,t)| \not= \alpha\\
\ge 0 \,\,\mbox{ if }\,\, u^*(x,t) = \rho_{min}\,\, \mbox{ or } \,\,(p^*(x,t) = -\alpha \,\,\mbox{ and }\,\, u^*(x,t) = 0)\\
\le 0 \,\,\mbox{ if }\,\, u^*(x,t) = \rho_{max}\,\, \mbox{ or }\,\, (p^*(x,t) = \alpha \,\,\mbox{ and } \,\,u^*(x,t) = 0)
\end{array}
\right.\,,
\end{equation}
as well as an analogous condition for $v_\Gamma$.
\Brem Let us compare the first condition in \eqref{pointwise} with the situation in the differentiable control problem without sparsity terms obtained for $\alpha = \aG=0$.
Then this condition boils down to the requirement that \,\,$v(x,t) = 0 \mbox{ if } \ |p^*(x,t) + \nu u^*(x,t)|> 0$,
or, since $\alpha=0$,
\begin{equation}  \label{simplecone}
v(x,t) = 0 \mbox{ if } \ |p^*(x,t) + \alpha \lambda^*(x,t) + \nu u^*(x,t)|> 0. 
\end{equation}
An analogous condition results for $\,\vG$.

One might be tempted to define the critical cone using \eqref{simplecone} and its counterpart for $\,\vG\,$ 
also in the case $\alpha > 0,\,\,\,\aG>0$.
This, however, is not a good idea, because it leads to a critical cone that is larger than needed, in general. 
As an example, we mention the particular case when the control
$\bus = \mathbf{0}$ satisfies the first-order necessary optimality conditions and when $\,\,|p^*| < \alpha\,\,$
and \,\,$|\psG|<\alpha_\Gamma$\,\, hold a.e. in $Q$ and $\Sigma$, respectively. 
Then the upper relation of \eqref{pointwise}, and its counterpart for $\vG$, lead to $C_{\bus}=\{\mathbf{0}\}$,
the smallest possible critical cone.

However, thanks to $u^*=0$, the variational inequality \eqref{varineq_full} implies that $p^* + \alpha \lambda^* + \nu u^*=0$ a.e. in $Q$, 
and hence the condition \,$|p^*(x,t) + \alpha \lambda^*(x,t) + \nu u^*(x,t)|> 0\,$ can only be satisfied on a set of measure zero. Moreover,
also the sign conditions \eqref{sign} do not restrict the critical cone, 
and therefore the largest possible critical cone $C_{\bus} = L^2(0,T;\CH)$ would be obtained, 
provided that analogous conditions hold for $u_\Gamma^*$ and $p^*_\Gamma$ on $\Sigma$.

In this example, the quadratic growth condition \eqref{growth} below is valid for the choice \eqref{critcone} as critical cone even 
without assuming the coercivity condition 
\eqref{coerc} below (here the so-called first-order sufficient conditions apply), while the use of a cone based on 
\eqref{simplecone} leads to postulating \eqref{coerc} on the whole space $L^2(0,T;\CH)$ for the quadratic growth 
condition to be valid. This shows that the choice of \eqref{critcone} as critical cone is essentially better than 
of one based on 
\eqref{simplecone}. 
\Erem

At this point, we give an explicit expression for $\,D^2\widehat J(\bu)[\bv,\mathbf{w}]\,$ for arbitrary $\,\bu=(u,\uG),
\bv=(v,\vG), \mathbf{w}=(w,w_\Gamma)\in\CU$. Arguing as in 
\cite[Sect.~5.7]{Fredibuch}, one obtains with $\,(\phi,\phi_\Gamma)=D\CS(\bu)[\bv]$\, and $\,(\psi,\psi_\Gamma)=D\CS(\bu)[\mathbf{w}]\,$
that
\begin{align}
\label{walter1}
&D^2\widehat J(\bu)[\bv,\mathbf{w}]=\iint_Q \bigl( \beta_1-p\, f^{(3)}(y)\bigr)\,\phi \psi\,+ \iint_\Sigma \bigl(\beta_2-p_\Gamma f_\Gamma^{(3)}(y_\Gamma)\bigr)
\,\phi_\Gamma \psi_\Gamma \nonumber\\
&\quad +\,\beta_3\iO \phi(T) \psi(T)\,+\,\beta_3\int_\Gamma \phi_\Gamma(T) \psi_\Gamma(T)\,+\,\nu\,\iint_Q  vw\,+\,\nu_\Gamma\,
\iint_\Sigma \vG w_\Gamma\,,
\end{align}  
where $\,(y,y_\Gamma)\,$ and $\,(p,p_\Gamma)\,$ are the state and the adjoint state associated with $\bu$. We claim that
\begin{align}
\label{walter2}
\left|D^2\widehat J(\bu)[\bv,\mathbf{w}]\right|\,\le\,\widehat C\,\|\bv\|_{L^2(0,T;\CH)}\,\|\mathbf{w}\|_{L^2(0,T;\CH)}\,,
\end{align}
where the constant $\widehat C>0$ is independent of $\,\bu,\bv,\mathbf{w}\in \Uad$. To prove the validity of \eqref{walter2}, 
we estimate the only critical
term 
$$
I:=\,-\iint_Q p f^{(3)}(y)\phi \psi \,-\iint_\Sigma p_\Gamma f_\Gamma^{(3)}(y_\Gamma) \phi_\Gamma \psi_\Gamma\,.
$$
To this end, recall that $(p,p_\Gamma)\in {\cal Y}\,$ by Theorem 4.3 and the global bound \eqref{ssbound2}.
Then, using H\"older's inequality, and the continuous
embeddings $\,V\subset L^4(\Omega)\,$ and $\,V_\Gamma\subset L^4(\Gamma)$, we obtain that
\begin{align*}
|I|\,&\le\,K_1\int_0^T\left( \|p\|_{L^2(\Omega)}\,\|\phi\|_{L^4(\Omega)}\,\|\psi\|_{L^4(\Omega)}
\,+\,\|p_\Gamma\|_{L^2(\Gamma)}\,\|\phi_\Gamma\|_{L^4(\Gamma)}\,\|\psi_\Gamma\|_{L^4(\Gamma)}\right)dt\\[2mm]
&\le\,C\left(\|\phi\|_{C^0([0,T];V)}\,\|\psi\|_{C^0([0,T];V)}\,+\,\|\phi_\Gamma\|_{C^0([0,T];V_\Gamma)}\,
\|\psi_\Gamma\|_{C^0([0,T];V_\Gamma)} \right)\\[2mm]
&\le\,C\,\|(\phi,\phi_\Gamma)\|_\CY \,\|(\psi,\psi_\Gamma)\|_\CY
\,\,\le\,\,C\,\|\bv\|_{L^2(0,T;\CH)}\,\|\mathbf{w}\|_{L^2(0,T;\CH)}\,,
\end{align*}
which proves the claim. This result shows that, for all $\bu \in \Uad$,  the functional $D^2\widehat J(\bu)$
can be continuously extended to a continuous bilinear functional on $L^2(0,T;\mathcal{H})^2$. This extension,
which will still be denoted by $D^2\widehat J(\bu)$, will be frequently used in the following.

We will rely on the following coercivity condition:

\begin{equation} \label{coerc}
D^2\widehat J(\bus)[\bv,\bv] >0 \quad \forall \, \bv \in C_\bus \setminus \{\mathbf{0}\}\,.
\end{equation}

\noindent Condition \eqref{coerc} is a direct extension of associated conditions that are standard in finite-dimensional 
nonlinear optimization. In the optimal control of partial differential equation, it was first used in \cite{casas_troeltzsch2012}.
As in  \cite[Thm 3.3]{casas_ryll_troeltzsch2015} or \cite{casas_troeltzsch2012}, it can be shown that \eqref{coerc} is equivalent to the existence of a 
constant $\delta > 0$ such that $D^2\widehat J(\bus)[\bv,\bv] \ge \delta \, \|\bv\|_{L^2(0,T;\CH)}^2$ for all $\bv \in C_{\bus}$.

We have the following result.
\Bthm \,\,{\rm (Second-order sufficient condition)} \,\,Suppose that {\bf (A1)}--{\bf (A6)} are fulfilled and that
$\,\rhomin<0<\rhomax\,$ and $\,\Rmin<0<\Rmax$. Moreover, let $\bus=(u^*,\usG) \in \Uad$, together with the associated state 
$(y^*,y_\Gamma^*)=\CS(\bus)$ and adjoint state $(p^*,p_\Gamma^*)$, fulfill the first-order necessary optimality conditions 
of Theorem 4.5. If, in addition,  $\bus$ satisfies the coercivity condition \eqref{coerc}, 
then 
there exist $\varepsilon > 0$ and $\sigma > 0$ such that the quadratic growth condition
\begin{equation} \label{growth}
\wJ(\bu) \ge \wJ(\bus) + \sigma \, \|\bu-\bus\|^2_{L^2(0,T;{\cal H})} 
\end{equation}
holds for all $\bu \in \Uad$ with $ \|\bu-\bus\|_{L^2(0,T;{\cal H})}  < \varepsilon$. Consequently,
$\bus$ is a locally optimal control in the sense of $L^2(0,T;{\cal H})$.
\Ethm
\begin{proof} The proof follows the one of \cite[Thm.~3.4]{casas_ryll_troeltzsch2015}. We remark that in 
\cite{casas_ryll_troeltzsch2015} the second-order differentiability of the objective functional in some $L^p$-space with $p < \infty$ 
was used, which we do not have in our situation. However, as E. Casas pointed out to us in a private communication,
this argument is not needed. 

We argue by contradiction, assuming that the claim of the theorem is not true. Then there exists a sequence of controls 
$\{\bu_k\}\subset \Uad$ such that, for all $k\in\enne$, 
\begin{equation}
\label{contrary}
\|\bu_k-\bus\|_{L^2(0,T;\CH)} < \frac{1}{k} \quad \mbox{ while } \quad \widehat{\mathcal{J}}(\bu_k)<  \widehat{\mathcal{J}}(\bus) 
+ \frac{1}{2k} \|\bu_k-\bus\|_{L^2(0,T;\CH)}^2\,.
\end{equation}
Noting that $\bu_k\not=\bus$ for all $k\in\enne$, we define 
$$
r_k = \|\bu_k-\bus\|_{L^2(0,T;\CH)}
  \quad \mbox{ and } \quad \bv_k = \frac{1}{r_k}(\bu_k-\bus)\,.
$$                                               
Then $\|\bv_k\|_{L^2(0,T;\CH)}=1$ and, possibly after selecting a subsequence, we can assume that 
\[
\bv_k \to \bv \, \mbox{ weakly in }\,  L^2(0,T;\CH)
\]
for some $\bv\in L^2(0,T;\CH)$. As in \cite{casas_ryll_troeltzsch2015}, the proof is split into three parts. 

(i) $\bv \in C_{\bus}$: Obviously, each $\bv_k$ obeys the sign conditions \eqref{sign} and thus belongs to $C(\bus)$. 
Since $C(\bus)$ is convex and closed in $L^2(0,T;\CH)$, it follows that $\bv\in C(\bus)$. We now claim that 
\Beq
\label{Paulchen}
D\widehat{J}(\bus)[\bv] + j'(\bus,\bv) = 0.
\end{equation}
Notice that by Remark 4.4 the expression 
$\,D\widehat{J}(\bus)[\bv]\,$ is well defined. For every $r \in (0,1)$ and all $\bv=(v,\vG),\, \bu=(u,\uG) \in L^2(0,T;\CH)$, 
we infer from the convexity of $j$ that
\begin{align}
j(\bv)-j(\bu) &\ge \frac{j(\bu + r (\bv-\bu))-j(\bu)}{r} \ge j'(\bu,\bv-\bu)\nonumber \\
&= \max_{(\alpha \lambda,\alpha_\Gamma\lambda_\Gamma)\in \partial j(\bu)}\Big(
\iint_Q\alpha \lambda (v-u) + \iint_\Sigma \alpha_\Gamma \lambda_\Gamma (v_\Gamma-u_\Gamma)\Big).
\label{directionalder}
\end{align}
This inequality yields, with $\bu_k=(u_k,u_{k_\Gamma})$,
\begin{align}
&D\widehat{J}(\bus)[\bv] +  j'(\bus,\bv)\ge D\widehat{J}(\bus)[\bv] + \iint_Q\alpha \lambda^* v + 
\iint_\Sigma \alpha_\Gamma \lambda_\Gamma^* v_\Gamma \nonumber\\
&= \iint_Q (p^* + \nu u^*) v +   \iint_\Sigma (p_\Gamma^* + \nu u_\Gamma^*) v_\Gamma +\iint_Q\alpha 
\lambda^* v + \iint_\Sigma \alpha_\Gamma \lambda_\Gamma^* v_\Gamma \nonumber\\
&= \lim_{k \to \infty} \frac{1}{r_k}  \Big(\iint_Q (p^* + \nu u^* + \alpha \lambda^*) (u_k-u^*)+\iint_\Sigma (p_\Gamma^* + \nu u_\Gamma^* + \alpha \lambda_\Gamma^*) ({u_k}_\Gamma-u_\Gamma^*)\Big)\nonumber\\
&\ge 0\,,
\end{align}
by the variational inequality \eqref{varineq_full}. Next, we prove the converse inequality. By \eqref{contrary}, we have
\[
\widehat{J}(\bu_k)-\widehat{J}(\bus) + j(\bu_k)-j(\bus) < \frac{1}{2k} r_k^2\,,
\]
whence, owing to the mean value theorem, and since $\bu_k = \bus+ r_k \bv_k$, we get
\[
\widehat{J}(\bus) + r_k D\widehat{J}(\bus + \theta_k r_k \bv_k)[\bv_k] + j(\bus + r_k \bv_k)
< \widehat{J}(\bus) + j(\bus) +  \frac{1}{2k} r_k^2
\]
with some $0 < \theta_k < 1$. From \eqref{directionalder}, we obtain $j(\bus + r_k \bv_k)-j(\bus)\ge  j'(\bus,r_k \bv_k)$, and thus
\[
r_k D\widehat{J}(\bus + \theta_k r_k \bv_k)[\bv_k] + r_k j'(\bus,\bv_k) < \frac{ r_k^2}{2k}\,.
\]
We divide this inequality by $r_k$ and pass to the limit $k \to \infty$. Here,  we invoke Corollary \ref{coroll3} of the Appendix, and
we use that $j'(\bus,\bv_k)\to j'(\bus,\bv)$. We then obtain the desired converse inequality
\[
D\widehat{J}(\bus)[\bv] + j'(\bus,\bv) \le 0\,,
\]
which completes the proof of (i).

(ii) $\bv = {\bf 0}$: We again invoke \eqref{contrary}, now  performing a second-order Taylor expansion on the left-hand side,
\begin{align}
&\widehat{J}(\bus) + r_k D \widehat{J}(\bus)[\bv_k] + \frac{r_k^2}{2} D^2 \widehat{J}(\bus + \theta_k r_k \bv_k)[\bv_k,\bv_k]
+ j(\bus + r_kv_k)\nonumber \\
&<\widehat{J}(\bus) + j(\bus) +  \frac{ r_k^2}{2k}\,.\nonumber 
\end{align}
We subtract $\widehat{J}(\bus) + j(\bus)$ from both sides and use \eqref{directionalder} once more to find that
\begin{equation} \label{Konrad2}
r_k  \left(D\widehat{J}(\bus)[\bv_k]  + j'(\bus,\bv_k)\right)+ \frac{r_k^2}{2} D^2 \widehat{J}(\bus + \theta_k r_k \bv_k)[\bv_k,\bv_k]<
\frac{ r_k^2}{2k}\,.
\end{equation}
From the right-hand side of \eqref{directionalder}, and the variational inequality \eqref{varineq2}, it follows
\[
D\widehat{J}(\bus)[\bv_k]  + j'(\bus,\bv_k) \ge 0\,,
\]
and thus, by \eqref{Konrad2},
\begin{equation}\label{liminf1} 
D^2\widehat{J}(\bus + \theta_k r_k \bv_k)[\bv_k,\bv_k]< \frac{1}{k}\,.
\end{equation}
Passing to the limit $k \to \infty$, we  apply Lemma \ref{LA3} and deduce that \,$D^2\widehat{J}(\bus)[\bv,\bv] \le 0.$
Since we know that $\bv \in C_{\bus}$, the second-order condition \eqref{coerc} implies that \,$\bv = \mathbf{0}$. 

(iii) {\em Contradiction:} To finish the proof, we employ \eqref{walter1} to see that
\begin{align}
&D^2\widehat J(\bus)[\bv_k,\bv_k]=\iint_Q \bigl( \beta_1-p^* f^{(3)}(\ys)\bigr)\,\phi_k^2\,+ \iint_\Sigma \bigl(\beta_2-\psG f_\Gamma^{(3)}(\ysG)\bigr)
\,\phi_{k_\Gamma}^2\nonumber\\
&\quad +\,\beta_3\iO \phi_k(T)^2\,+\,\beta_3\int_\Gamma \phi_{k_\Gamma}(T)^2\,+\,\nu\,\iint_Q v_k^2\,+\,\nu_\Gamma\,
\iint_\Sigma v_{k_\Gamma}^2. \label{Konrad}
\end{align}
As shown in the previous step, $\bv = \mathbf{0}$, and therefore $\,\bv_k \to \mathbf{0}$ weakly in $L^2(0,T;\CH)$. By 
Lemma \ref{LA3}, the sum of the four integrals containing $\phi_k$ or $\phi_{k_\Gamma}$  tends to zero. 
On the other hand, we have $\,\|\bv_k\|_{L^2(0,T;\CH)} = 1\,$ for all $k\in\enne$, by construction. Hence, 
\begin{equation}
\label{auweia!}
\nu\,\iint_Q v^2_k\,+\,\nu_\Gamma\,\iint_\Sigma v_{k_\Gamma}^2 \,\ge\, \min \{\nu,\nu_\Gamma\} 
\,\left(\iint_Q v^2_k + \iint_\Sigma v_{k_\Gamma}^2 \right) 
= \min \{\nu,\nu_\Gamma\} > 0.
\end{equation}
It therefore follows from the weak sequential lower semicontinuity of the last two summands on the right-hand side of \eqref{Konrad}
that 
\begin{align*}
&\liminf_{k\to\infty} \,D^2\widehat J(\bus)[\bv_k,\bv_k] \,\ge \,\liminf_{k\to\infty}\Big(
\nu\,\iint_Q v^2_k\,+\,\nu_\Gamma\,\iint_\Sigma v_{k_\Gamma}^2\Big)\,\ge\,\min \{\nu,\nu_\Gamma\}\,>0\,.
\end{align*}
On the other hand, it is easily deduced from \eqref{liminf1} and \eqref{lip2} that 
\[
\liminf_{k \to \infty} D^2\widehat{J}(\bus)[\bv_k,\bv_k] \le 0\,,
\]
a contradiction. The assertion of the theorem is thus proved. 
\end{proof}

For the particular case $\alpha=\alpha_\Gamma$ without sparsity functional, Theorem 4.10 improves the second-order 
sufficient condition  \cite[Thm.~3.6]{CS}: indeed, our coercivity condition \eqref{coerc} is required on a smaller 
critical cone (compare \eqref{simplecone} with the condition \cite[(3.72)]{CS}), and we have local optimality in an 
$L^2$-neighborhood, hence in a larger set than in an $L^\infty$-neighborhood as in \cite{CS}.

\section{Appendix}
\setcounter{equation}{0}
In the following, we assume that {\bf (A1)}--{\bf (A6)} are fulfilled.
\Blem 
\label{LA1} 
Let $\{\bu_k\}\subset\Uad$ converge strongly in $L^2(0,T;\CH)$ to $\bus\in\Uad$. Then the sequence  $\{\by_k\}$ 
of associated states converges strongly in $\,\CY\,$ to $\,\bys$, and the sequence $\{\bp_k\}$ of associated 
adjoint states converges strongly in $\,\CY\,$ to $\bp^*$.
\Elem
\begin{proof} The strong convergence $\,\|\by_k-\bys\|_\CY\to 0\,$ follows directly from \cite[Lem.~2.4]{CS}. By the continuity
of the embedding $\,\CY\subset C^0([0,T];\CV),$ we then have \,$\|\by_k(T)-\bys(T)\|_\CV\to 0$. Moreover, since the states 
$\,\by_k=(y_k,y_{k_\Gamma})\,$ and $\,\bys=(\ys,\ysG)\,$ have to obey the separation property \eqref{separation}, we can easily infer 
from \eqref{ssbound2} and the continuous
embedding $\CV\subset \left(L^6(\Omega)\times L^6(\Gamma)\right)$, using the mean value theorem, that
\begin{equation}
\label{Hugo1}
\|f''(y_k)-f''(\ys)\|_{C^0([0,T];L^6(\Omega))}\,+\,\|f_\Gamma''(y_{k_\Gamma})-f_\Gamma''(\ysG)\|_{C^0([0,T];L^6(\Gamma))}
\to 0 \quad\mbox{as }\,k\to\infty.
\end{equation}
Next, we observe that the
adjoint states $\bp_k=(p_k,{p_k}_\Gamma)$ solve the system
\begin{align}
\nonumber
&-\dt p-\Delta p + f''(y_k)p =\beta_1(y_k-y_Q) \quad\mbox{a.e. in }\,Q,\\
&-\dt \pG -\delG\pG+\dn  p  + f_\Gamma''({y_k}_\Gamma)\pG=\beta_2({y_k}_\Gamma-y_\Sigma) 
\quad\mbox{and} \quad p_\Gamma=p_{|\Gamma} \quad\mbox{a.e. on }\,\Sigma,\nonumber \\
\nonumber
&p(T)=\beta_3(y_k(T)-y_{\Omega,T}) \quad\mbox{a.e. in }\,\Omega, \quad \pG(T)=\beta_3({y_k}_\Gamma(T)-y_{\Gamma,T}) \quad\mbox{a.e. on }\,\Gamma. \nonumber 
\end{align}
From \eqref{ssbound2} it follows that the sequences $\,\{\|f''(y_k)\|_{L^\infty(Q)}\}\,$ and 
$\,\{\|f_\Gamma''({y_k}_\Gamma)\|_{L^\infty(\Sigma)}\}\,$
are bounded. Arguing as in the proof of the bound \eqref{adj4} in Theorem 4.3, we obtain that
\begin{align}
\|\bp_k\|_{\mathcal{Y}} 
\le c\left(\|y_k-y_\Omega\|_{L^2(Q)} + \|{y_k}_\Gamma-y_\Sigma\|_{L^2(\Sigma)} 
+\|y_k(T)-y_{\Omega,T}\|_V + \|{y_k}_\Gamma(T)-y_{\Gamma,T}\|_{V_\Gamma}\right) \nonumber
\end{align}
for all $k\in\enne$. In view of the convergence results shown above, we thus can conclude that
\begin{equation} \label{boundpk}
\|\bp_k\|_{\mathcal{Y}} \le K \quad \mbox{for all $k\in\enne$, with some constant $\,K>0$}. 
\end{equation}
Now we subtract the adjoint equations for $\bp_k$ and $\bps$ and set
$\bz_k= (z_k,z_{k_\Gamma})=\bp_k - \bps$. After some rearrangement, we arrive at the system
\begin{align}
&-\dt z_k-\Delta z_k + f''(y^*)z_k =\beta_1(y_k-y^*) + [f''(y_k)-f''(y^*)]p_k \quad\mbox{a.e. in }\,Q,\label{z1}\\
&-\dt {z_k}_\Gamma -\delG{z_k}_\Gamma+\dn  {z_k}  + f_\Gamma''({y}_\Gamma^*){z_k}_\Gamma=\beta_2({y_k}_\Gamma-y_\Gamma^*) 
+ [f_\Gamma''({y_k}_\Gamma)-f_ \Gamma''(y_\Gamma^*)]{p_k}_{\Gamma} \nonumber\\
&\qquad\mbox{and }\, {z_k}_\Gamma={z_k}_{|\Gamma}\quad\mbox{a.e. on }\,\Sigma, \label{z2} \\
& z_k(T)=\beta_3(y_k(T)-y^*(T))\quad\mbox{a.e. in }\,\Omega, \quad {z_k}_\Gamma(T)=\beta_3({y_k}_\Gamma(T)-y_{\Gamma}^*(T))
\quad\mbox{a.e. on }\,\Gamma\,. \label{z4} 
\end{align}
Again, we apply  Theorem 4.3 to estimate $\bz_k$ in terms of the norms of the right-hand sides. Now notice that
from \eqref{Hugo1} and \eqref{boundpk} it readily follows that the right-hand side of \eqref{z1} converges to zero strongly in $L^2(Q)$. 
Analogously, the right-hand side of \eqref{z2}
tends to zero strongly in $L^2(\Sigma)$. Therefore, and since $\,\|\by_k(T)-\bys(T)\|_\CV\to 0$, we can infer from Theorem 4.3 that
 $\,\|\bz_k\|_\CY\to 0\,$ as $k\to\infty$. The assertion is thus proved.
\end{proof}

\begin{corollary} \label{coroll3} Let $\{\bu_k\}\subset\Uad$ converge strongly in $L^2(0,T;\CH)$ to $\bus\in\Uad$, and let
$\{\bv_k\}$ converge weakly to $\bv$ in $L^2(0,T;\CH)$. Then 
\begin{equation}
\label{Hugo2}
\lim_{k\to\infty} \,D\widehat{J}(\bu_k)[\bv_k] \,=\, D\widehat{J}(\bus)[\bv]\,.
\end{equation}
\end{corollary}
\begin{proof} We have, with $\bv_k=(v_k,v_{k_\Gamma})$,
\[
D\widehat{J}(\bu_k)[\bv_k] = \iint_Q(p_k + \nu u_k)v_k + \iint_\Sigma ({p_k}_\Gamma + \nu {u_k}_\Gamma){v_k}_\Gamma.
\]
Owing to Lemma \ref{LA1}, we have, in particular, that  $\,\{\bp_k + \nu \bu_k\}\,$ converges to $\,\bps+\nu\bus\,$ strongly in
$\,L^2(0,T;\CH)$, whence the assertion immediately follows.
\end{proof}
\begin{lemma} \label{LA3} Let $\{\bu_k\}$ and $\{\bv_k\}$ satisfy the conditions of Corollary \ref{coroll3}, and
assume that $\,\nu=\nu_\Gamma=0$. Then
\begin{equation} \label{liminf0}
\lim_{k \to \infty} D^2\widehat{J}(\bu_{k})[\bv_k,\bv_k] = D^2\widehat{J}(\bus)[\bv,\bv].
\end{equation} 
\end{lemma}
\begin{proof}
Let $\bv_k=(v_k,v_{k_\Gamma})$, $\bv=(v,\vG)$, $\,(\varphi_k,\varphi_{k_\Gamma})=D\CS(\bu_k)[\bv_k]$, and $(\varphi,\varphi_\Gamma)
=D\CS(\bus)[\bv]$. Since $\nu=\nu_\Gamma=0$, we infer from \eqref{walter1} that
\begin{align}
&D^2\widehat{J}(\bu_{k})[\bv_{k},\bv_{k}]=\iint_Q (\beta_1-p_k f^{(3)}(y_k))\varphi_k^2+\iint_\Sigma (\beta_2-{p_k}_\Gamma 
f_\Gamma^{(3)}({y_k}_\Gamma)){\varphi_k^2}_\Gamma \nonumber\\ 
&\qquad+ \beta_3 \int_\Omega \varphi_k^2(T) + \beta_3 \int_\Gamma  {\varphi_k^2}_\Gamma(T)
\,=\,\sum_{i=1}^4 I_{i,k},
\nonumber 
\end{align}
with obvious notation. At first, notice that
$$
(\varphi_k,{\varphi_k}_\Gamma)-(\varphi,\varphi_\Gamma) = \left(D\mathcal{S}(\bu_k)-D\mathcal{S}(\bus)\right)[ \bv_k]\, +\,D\CS(\bus)[\bv_k-\bv]\,.
$$
By virtue of \eqref{lip1} (recall Remark 3.1 in this regard) and the boundedness of $\{\bv_k\}$ in $L^2(0,T;\CH)$, the 
first summand on the right converges strongly to zero in $\,\CY$.
The second converges to zero weakly in $\,\CY$ and, thanks to the compactness of the embedding $\,\CY\subset
C^0([0,T];L^p(\Omega)\times L^p(\Gamma))$\, for $1\le p<6$ (see, e.g., \cite[Sect.~8, Cor.~4]{Simon1987}), strongly in
$\,C^0([0,T];L^5(\Omega)\times L^5(\Gamma))$. In conclusion,
\begin{equation}
\label{Hugo3}
(\varphi_k,\varphi_{k_\Gamma})\to (\varphi,\varphi_\Gamma)\quad\mbox{strongly in }\,C^0([0,T];L^5(\Omega)\times L^5(\Gamma))\,.
\end{equation}
In particular,
\begin{equation} \label{Hugo4}
\lim_{k\to \infty} \left(I_{3,k}+I_{4,k}\right)= \beta_3 \int_\Omega \varphi^2(T) + \beta_3 \int_\Gamma  \varphi^2_\Gamma(T)\,.
\end{equation}
Moreover, similarly as in \eqref{Hugo1}, we have, as $k\to\infty$,
\begin{equation}
\label{Hugo5}
\|f^{(3)}(y_k)-f^{(3)}(\ys)\|_{C^0([0,T];L^6(\Omega))} + \|f^{(3)}_\Gamma(y_{k_\Gamma})-f^{(3)}_\Gamma(\ysG)\|_{C^0([0,T];L^6(\Gamma))}
\to 0\,,
\end{equation}
and we know already from Lemma 5.1 that $\,\bp_k\to\bps\,$ strongly in $C^0([0,T];L^6(\Omega)\times L^6(\Gamma))$.
Combining this with \eqref{Hugo3} and \eqref{Hugo5}, and invoking H\"older's inequality appropriately, we easily verify that
\begin{equation} \label{conv1}
\lim_{k \to \infty}(I_{1,k}+I_{2,k}) = \iint_Q (\beta_1-p^* f^{(3)}(y^*)) \varphi^2\,+\iint_\Sigma (\beta_2-\psG f_\Gamma^{(3)}(\ysG))
\varphi_\Gamma^2\,.
\end{equation} 
From \eqref{Hugo4} and \eqref{conv1}, the assertion follows.
\end{proof}
\noindent
{\bf Acknowledgment.} The authors thank Prof. Eduardo Casas (Santander, Spain) for encouraging them to use the critical cone
$C_{\bus}$ in place of an extended cone as in \cite{CS}. His idea paved the way for adapting the results of 
\cite{casas_ryll_troeltzsch2015} to this lower regularity case.


\End{document}